\documentclass[final,leqno]{siamltex}
\pdfoutput=1
\usepackage{amsmath, amsfonts}
\usepackage{caption}
\usepackage{subcaption}
\usepackage{graphicx,color}
\usepackage{multirow}

\newtheorem{remark}{Remark}

\newcommand{\bvec}[1]{\mathbf{#1}}

\renewcommand{\Re}{\mathrm{Re}~}
\renewcommand{\Im}{\mathrm{Im}~}
\newcommand{\Tr}{\mathrm{Tr}}

\newcommand{\vr}{\bvec{r}}

\newcommand{\ud}{\,\mathrm{d}}

\newcommand{\Vxc}{\hat{V}_{\mathrm{xc}}}
\newcommand{\Vion}{\hat{V}_{\mathrm{ion}}}
\newcommand{\abs}[1]{\lvert#1\rvert}

\newcommand{\ie}{\textit{i.e.}}
\newcommand{\eg}{\textit{e.g.}}
\newcommand{\Or}{\mathcal{O}}

\newcommand{\lmin}{\lambda_{\min}}
\newcommand{\lmax}{\lambda_{\max}}
\newcommand{\Ran}{\text{Ran}}

\newcommand{\response}[1]{{#1}}

\title{Pole expansion for solving a type of parametrized linear systems
in electronic structure calculations}

\author{
Anil Damle \thanks{Institute for Computational and Mathematical Engineering,
Stanford University, Stanford, CA 94305. Email: damle@stanford.edu} \and
Lin Lin \thanks{
Computational Research Division, Lawrence Berkeley National
Laboratory, Berkeley, CA 94720. Email: linlin@lbl.gov} 
\and Lexing Ying \thanks{Department of Mathematics and 
Institute for Computational and Mathematical Engineering,
Stanford University, Stanford, CA 94305. Email: lexing@math.stanford.edu}
}

\begin{document}

\maketitle

\begin{abstract}
	We present a new method for solving parametrized linear systems.
	Under certain assumptions on the parametrization, 
	solutions to the linear systems for all parameters can be accurately
	approximated by linear combinations of solutions to linear systems
	for a small set of fixed parameters.  Combined with either direct
	solvers or preconditioned iterative solvers for each linear system
	with a fixed parameter, the method is particularly suitable for situations when
	solutions to a large number of distinct parameters or a large number of right
	hand sides are required.  The method is also simple to parallelize.  We
	demonstrate the applicability of the method to the calculation of the
	response functions in electronic structure theory.  We demonstrate the
	numerical performance of the method using a benzene molecule and a DNA
	molecule.
\end{abstract}


\pagestyle{myheadings}
\thispagestyle{plain}
\markboth{A. DAMLE, L. LIN AND L. YING}{POLE EXPANSION FOR
PARAMETRIZED LINEAR SYSTEMS}


\section{Introduction}\label{sec:intro}

Consider the linear system 
\begin{equation}
	(H - zS) u = b,
	\label{eqn:multishift}
\end{equation}
where $H$ and $S$ are $N\times N$ Hermitian matrices, \response{$S$ is positive definite}, $u$ and $b$ are vectors of
length $N$, and $z\in \mathbb{C}$ with $\Re z\le 0$. 

Often~\eqref{eqn:multishift} needs to be solved a large number of
times, both for a large number of distinct parameters $z$, and for
multiple right hand sides $b$.  This type of calculation arises in a
number of applications, such as time-independent density functional
perturbation theory
(DFPT)~\cite{BaroniGiannozziTesta1987,GonzeAllanTeter1992,BaroniGironcoliDalEtAl2001},
many body perturbation theory using the GW
method~\cite{Hedin1965,AryasetiawanGunnarsson1998,FriedrichSchindlmayr2006,UmariStenuitBaroni2010,PingRoccaGalli2013,GiustinoCohenLouie2010},
and the random phase approximation (RPA) of the electron correlation
energy~\cite{LangrethPerdew1975,LangrethPerdew1977,Furche2001,NguyenGironcoli2009}.
The connection between~\eqref{eqn:multishift} and these
applications will be illustrated in Section~\ref{sec:connection}.

\subsection{\response{Previous work}}
From an algorithmic point of view, solving linear systems with multiple
shifts has been widely explored in the
literature~\cite{Frommer2003,FrommerGlassner1998,SimonciniSzyld2007,DattaSaad1991,Meerbergen2003,GallivanGrimmeVan1996,BaiFreund2001,FeldmannFreund1995,SABAKI2012}.
These methods are often based on the Lanczos method.  The basic idea here is
that when $S=I$ the Krylov subspace constructed in the Lanczos procedure
is invariant to the shifts $z$. Therefore the same set of Lanczos
vectors can be used simultaneously to solve~\eqref{eqn:multishift}
with multiple $z$. When $S$ is not the identity matrix, $S$ should be factorized,
such as by the Cholesky factorization, and \response{triangular solves with the Cholesky factors must be
computed} during every iteration, which can significantly increase the
computational cost.  

In its simplest form the Lanczos method as described in Section
\ref{sec:lanczos} precludes the use of a preconditioner. However, there
are strategies, some of which are discussed in the aforementioned
references, for using a small number of preconditioners and perhaps multiple
Krylov spaces to find solutions for all of the desired parameters. Even
when using preconditioners it may be unclear how to pick which Krylov
spaces to use, or what the shifts for the preconditioners should be, see
\eg, \cite{Meerbergen2003}. Another specific class of methods are the so called recycled Krylov methods, see, \eg, \cite{Parks2006}. In such methods, the work done to solve one problem, potentially with a preconditioner, is reused to try and accelerate the convergence of related problems. 

Another class of popular methods, especially in the context of model
reduction, are based on the Pad\'{e} approximation, which does not target the
computation of all the entries of $u(z)=(H-zS)^{-1} b$ as in~\eqref{eqn:multishift}, but rather a linear functional of $u(z)$ in the form
of $f(z)=l^{T}u(z)$.  $f(z)$ is a scalar function of $z$,
which can be stably
expanded in the Pad\'{e} approximation via the Lanczos
procedure~\cite{FeldmannFreund1995}.  However, it is not obvious how to
obtain a Pad\'{e} approximation for all entries of $u(z)$ directly, and with an
approximation that has a uniformly bounded error for all $\Re z\le 0$.

\subsection{\response{Contribution}}
In this paper we develop a new approach for solving para\-metrized
linear systems using a pole expansion. The pole expansion used here is a rational
approximation to all entries of $u(z)$ simultaneously. The main idea behind this
pole expansion comes directly from the work of Hale, Higham and
Trefethen~\cite{HHT}, which finds a nearly optimal quadrature rule of
the contour integral representation for $u(z)$ in the complex plane.
The idea of the work in~\cite{HHT} has also been adapted in the context of
approximating the spectral projection operator and the Fermi-Dirac
operator (a ``smeared'' spectral projection operator) for ground state electronic structure
calculation~\cite{LinLuYingE2009}. 

Such a scheme directly expresses $u(z),$ and thus the solution to \eqref{eqn:multishift} for any $z,$
as the linear combination of solutions to a small set of linear systems, each of which has a fixed parameter.
Each linear system with a fixed parameter may be solved either with a preconditioned iterative
method or a direct method. Furthermore, the construction of the pole expansion depends only on the largest and smallest positive generalized eigenvalues of the matrix pencil and the number of poles. The number of poles used in the expansion yields control over the accuracy of the scheme. In fact, the pole expansion converges exponentially
with respect to the number of approximating terms across the range
$\Re z \le 0.$  Finally, the pole expansion allows the same
treatment for general $S$ as in the case of $S=I,$ without directly
using its Cholesky factor  or $S^{-1}$ at every iteration.  As a result,
the pole expansion addresses the disadvantages of both the Lanczos method
and the Pad\'{e} approximation.

\subsection{Notation}
\begin{figure}
  \centering
  \includegraphics[width = 0.6\textwidth]{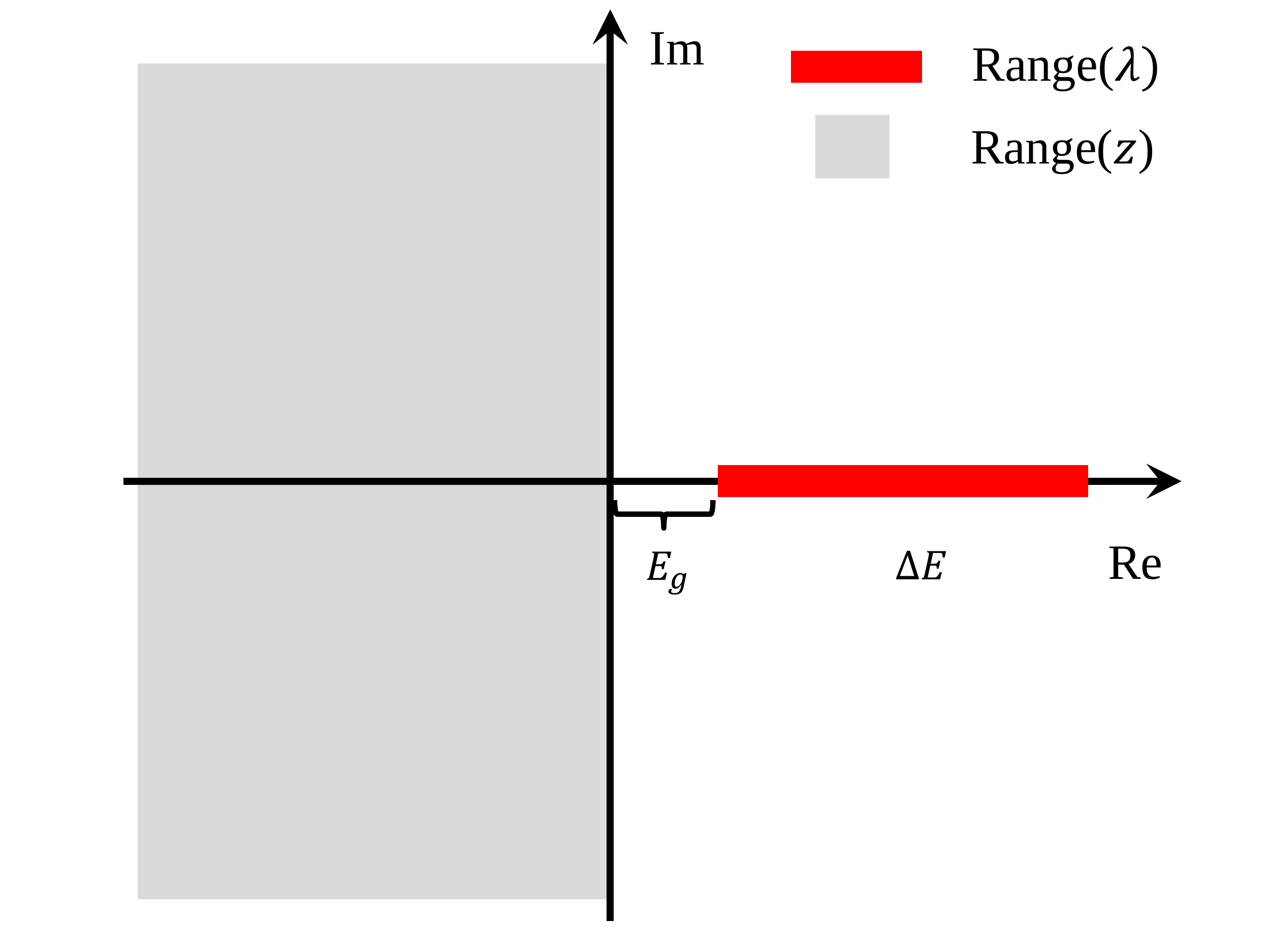} 
  \caption{A schematic view of the range of the spectrum of the
	$(H,S)$ \response{pencil} (thick red line on the positive real axis),
	and the range of the parameter $z$ (light gray area), separated by a
	positive distance $E_{g}$.}
  \label{fig:pole1}  
\end{figure}

In this paper we use the following notation. Let
$\{\lambda_{i}\},\{\psi_{i}\}$ be the generalized eigenvalues and eigenvectors of
the matrix \response{pencil} $(H,S)$ which satisfy
\begin{equation}
\label{eqn:GE}
	H\psi_{i} = S \psi_{i} \lambda_{i},\quad i=1,\ldots,N,
\end{equation}
with $\{\lambda_{i}\}$ labeled in non-increasing order. Such an
eigen decomposition is only used for deriving the pole expansion, and is
not performed in practical calculations.
We initially assume that all the eigenvalues $\lambda_{i}$ are positive, and denote by
$E_{g}\equiv \lambda_N>0$. However, we later relax this assumption and allow for a small number of negative eigenvalues. We also denote by $\Delta E$ the spectrum
width of the $(H,S)$ \response{pencil} \ie\  $\Delta E = \lambda_{1}-\lambda_{N}$.
We refer to Fig.~\ref{fig:pole1} for an illustration of the relative position
of the spectrum and the range of the parameter $z$ in the complex plane.

Furthermore, $A^*$ and $A^{T}$
denote the conjugate transpose and  the transpose of a matrix $A$,
respectively. We let $R$ be an upper triangular matrix that denotes the
Cholesky factorization of $S,$ \ie, 
\begin{equation*} 
	S = R^* R.
\end{equation*}
We use the notation
\begin{equation*}
\mathcal{K}_{n}(A,b)  = \text{span}\left(b,Ab,\ldots,A^{n-1}b \right)
\end{equation*}
to denote the $n^{th}$ Krylov subspace associated with the matrix $A$ and
the vector $b.$ Finally, we let $N_z$ be the number of distinct complex
shifts $z$ for which we are interested in solving~\eqref{eqn:multishift}.

\subsection{Outline}
The rest of the manuscript is organized as follows.  We first discuss
the standard Lanczos method for solving linear systems with multiple
shifts in Section~\ref{sec:lanczos}.  We introduce the pole expansion method
for solving~\eqref{eqn:multishift}, and  analyze
the accuracy and complexity of the approach in Section~\ref{sec:pole}. In
general the matrix \response{pencil} $(H,S)$ may not always have all positive eigenvalues,
and the case where there are positive and negative eigenvalues is discussed in
Section~\ref{sec:indefinite}.
The connection between~\eqref{eqn:multishift} and electronic
structure calculations is given in Section~\ref{sec:connection}.
The numerical results with applications
to density functional perturbation theory calculations are given in
Section~\ref{sec:numer}.
We conclude in Section~\ref{sec:conclusion}.


\section{Lanczos method for parametrized linear systems}\label{sec:lanczos}
We briefly describe \response{a basic variant of the} Lanczos method for solving parametrized linear systems of the form  \response{\eqref{eqn:multishift}}
where the matrix \response{pencil} $(H,S)$ satisfies the conditions given in the
introduction.
Using the Cholesky factorization of $S$ we may transform~\eqref{eqn:multishift} in a manner such that we instead solve 
\begin{equation}
\label{eqn:chol_eq}
(R^{-*}HR^{-1}-zI)\tilde{u} = \tilde{b},
\end{equation}
where $u = R^{-1}\tilde{u}$ and $\tilde{b} = R^{-*}b.$ Since $R^{-*}HR^{-1}$ is Hermitian positive definite we now briefly describe the Lanczos method for parametrized systems under the assumption that our equation is of the form
\begin{equation*}
(A-zI)x = \tilde{b},
\end{equation*}
for some Hermitian positive definite matrix $A,$ right hand side $\tilde{b}$ and complex shift $z$ such that $\Re z \leq 0.$ 

Given these systems, we note that
\begin{equation*}
\mathcal{K}_k(A,\tilde{b}) = \mathcal{K}_k(A-zI,\tilde{b})
\end{equation*}
for any complex scalar $z.$ Based on this observation, the well known Lanczos process, see, \eg, \cite{GVL} for details \response{and \cite{Lanczos1950} for a historical perspective}, may be slightly modified such that at each iteration approximate solutions, each of which satisfy the Conjugate Gradient (CG) error criteria, are produced for each shift. \response{However, in this basic formulation using the invariance of the Krylov subspaces precludes the use of preconditioners for the problems.}

\begin{remark}
 We note that the restriction $\Re z \leq 0$ may be relaxed if we
 instead use a method analogous to MINRES to simultaneously solve the
 shifted systems, for details of MINRES see, \eg, \cite{PASA1975}. Such
 an algorithm is only slightly more expensive than the CG style
 algorithm given here and scales asymptotically in the same manner. 
\end{remark}

We observe that the computational cost of using this method breaks down
into two distinct components. First, there is the cost of the Lanczos
procedure which only has to be computed once regardless of the number of
shifts. The dominant cost of this procedure is a single matrix vector
multiplication at each iteration. In the case where $S\neq I$ the method
actually requires a single application of $R^{-*}HR^{-1},$ which is
accomplished via two triangular solves and one matrix vector product
with $H.$ Second, for each shift $z$ a $k\times k$ tridiagonal system, denoted $T_k,$
must be solved. However, this may be done very cheaply by maintaining a
$LDL^*$ factorization of $T_k$ for each shift (see, \eg, \cite{GVL}
for a detailed description of the case without a shift). This means that the 
computational costs at each iteration for solving the set of
sub-problems scales linearly with respect to both the
number of shifts and the problem size. Thus, the dominant factor in the
computation is often the cost of the Lanczos process. 

Furthermore, properties of the Lanczos process and the structure of the $LDL^*$
factorization imply that if we neglect the cost of storing $A$ the memory
costs of the algorithm scale linearly in both the number of shifts and
the problem size. At the core of this memory scaling is the fact that
not all of the Lanczos vectors must be stored to update the solution, see, \eg, \cite{GVL}. In fact, for
large problems that take many iterations storing all of the vectors
would be infeasible. However, this does impact the use of the algorithm
for multiple shifts. Specifically, the set of shifts for which solutions are desired must be decided
upon before the algorithm is run so that all of the necessary
factorizations of $T_k$ may be built and updated at each iteration. If a
solution is needed for a new shift the algorithm would have to be run
again from the start unless all of the Lanczos vectors were stored.

\section{Pole expansion for parametrized linear systems}\label{sec:pole}
We now describe a method based on the work of Hale, Higham and Trefethen \cite{HHT} to simultaneously solve systems of the form 
\begin{equation}
(H-z_lS)u_l = b
\label{eqn:pole_shift}
\end{equation}
\response{where the subscript $l$ explicitly denotes the set of $N_z$ distinct shifts. Here, the matrix pencil $(H,S)$ and the shifts $z_l$ satisfy the same properties as in the introduction.}

\subsection{Constructing a pole expansion}

First, we provide a brief overview of the method for computing functions
of matrices presented in \cite{HHT}. Given an analytic function $f$,
a Hermitian matrix $A,$ and a closed contour $\Gamma$ that
only encloses the analytic region of $f$ and winds once around the spectrum of
$A$ in a counter clockwise manner, then $f(A)$ may be represented via contour
integration as \begin{equation}
\label{eqn:contour}
f(A) = \frac{1}{2\pi i}\int_{\Gamma}f(\xi)(\xi I-A)^{-1}d\xi.
\end{equation}
The authors in \cite{HHT} provide a method via a sequence of conformal mappings to generate an efficient quadrature scheme for~\eqref{eqn:contour} based on the trapezoidal rule. Specifically, in \cite{HHT} a map from the region $S = [-K,K]\times [0,K']$ to the upper half plane of ${\Omega = \mathbb{C} \backslash ((-\infty,0]\cup [\lmin (A), \lmax (A)])}$ is constructed via 
\begin{equation*}
  \begin{split}
    z &= \sqrt{\lmin (A) \lmax (A)}\left(\frac{k^{-1}+u}{k^{-1}-u} \right),\\
    u &= sn(t) = sn(t|k),\\
    k &= \frac{\sqrt{\lmax (A) / \lmin (A)}-1}{\sqrt{\lmax (A) / \lmin (A)}+1},
  \end{split}
\end{equation*}
where $t\in S$ and $z$ is in the upper half plane of $\Omega.$ The constants $K$ and $K'$ are complete elliptic integrals and $sn(t)$ is one of the Jacobi elliptic functions. Application of the trapezoidal rule in $S$ using the $P$ equally spaced points 
\begin{equation*}
t_j = -K + \frac{iK'}{2}+2\frac{(j-\frac{1}{2})K}{P}, \;\; 1\leq j \leq P,
\end{equation*}
yields a quadrature rule for computing~\eqref{eqn:contour}. In \cite{HHT} the assumption made is that the only non-analytic region of $f$ lies on the negative real axis. 

Here, we instead consider the case where the non-analytic region of $f$
may be anywhere in the negative half plane. Therefore,
for our purposes a modification of the transform in \cite{HHT} must be
used.  Specifically, we use the construction of the quadrature presented
in Section 2.1 of \cite{LinLuYingE2009}, where now an additional
transform of the form $\xi = \sqrt{z}$ is used to get the quadrature
nodes. We note that there is a slight difference between the contour
used here and the one is \cite{LinLuYingE2009}. Because we are assuming
the matrix is Hermitian positive definite we only need to consider a
single branch of the square root function in defining the nodes $\xi_j,$
in this case the positive one. 

 The procedure outlined may be used to generate a $P$ term pole expansion for a function $f(A)$ denoted
\begin{equation}
\label{eqn:Pole_expand} 
f_P(A) \approx \sum_{k=1}^P w_k f(\xi_k)\left(\xi_k I - A\right)^{-1},
\end{equation} 
where $\xi_k$ and $w_k,$ $k=1,\ldots,P$ are the quadrature nodes and weights respectively. Both the nodes and the weights depend on $P,$ $\lmax (A),$ and $\lmin (A)$. In \cite{HHT} asymptotic results are given for the error $\|f(A)-f_P(A)\|.$ For a Hermitian positive definite matrix $A$ the asymptotic error in the expansion of the form~\eqref{eqn:Pole_expand}, see \cite{HHT} Theorem 2.1 and Section 2 of \cite{LinLuYingE2009}, behaves as  
\begin{equation}
\label{eqn:HHTerror}
\|f(A)-f_P(A)\| = \Or \left( e^{-C P / \log\left(\lambda_{max}(A)/\lambda_{min}(A)\right)} \right),
\end{equation}
where $C$ is a constant independent of $P,\lambda_{max}(A)$ and $\lambda_{min}(A).$ 

\begin{remark}
We make special note of the fact that in the framework where we
have all positive generalized eigenvalues the constant in the exponent
does not depend on \response{where the function $f$ has poles in the left half plane.} Specifically, the
initial use of the $z = \xi^2$ mapping maps the poles on the imaginary
axis to the negative real axis, see \cite{LinLuYingE2009} for details.
All of the poles \response{with strictly negative real part} get mapped to a distinct sheet
of the Riemann surface from the one that contains the positive real axis
and the poles that were initially on the imaginary axis. As noted
previously, the only sheet we need to consider building a contour on is
the one containing the spectrum of $A$, and thus $C$ does not depend on the \response{locations of the poles.}
\end{remark}

\begin{figure}
  \centering
  \includegraphics[width = 0.6\textwidth]{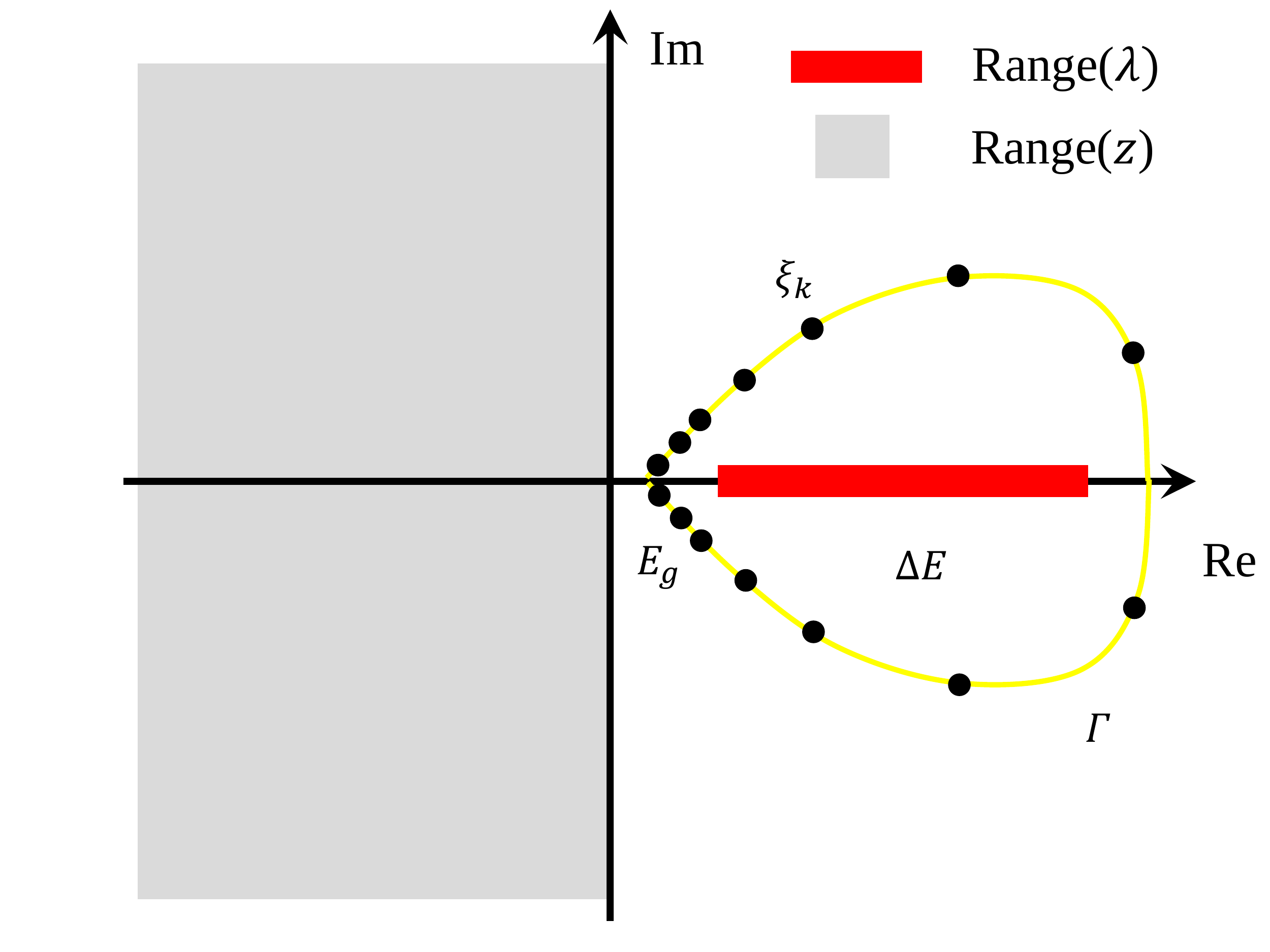} 
  \caption{A schematic view of the range of the spectrum of the
	$(H,S)$ \response{pencil} (thick red line on the positive real axis),
	and the range of the parameter $z$ (light gray area), separated by a
	positive distance $E_{g}$, together of the Cauchy contour (yellow
	line) surrounding the positive part of the spectrum and the
	discretized contour integration points (poles, black dots).}
  \label{fig:pole2}  
\end{figure}

Fig.~\ref{fig:pole2} illustrates the type of contours computed by this scheme and used throughout the remainder of the paper. To facilitate the computation of certain quantities necessary in the generation of these contours we used a Schwarz-Christoffel toolbox \cite{DriscollSC1,DriscollSC2}.  

\subsection{\response{Solving shifted linear systems}}
We now address the use of a pole expansion to solve problems of the form~\eqref{eqn:pole_shift}. 
Denote by $\Psi = [\psi_{1},\cdots,\psi_{N}]$, and $\Lambda =
\mathrm{diag}[\lambda_{1},\cdots,\lambda_{N}]$, from
the generalized eigenvalue problem ~\eqref{eqn:GE} we then have
\begin{equation*}
	\Psi^{*} H \Psi = \Lambda, \quad \Psi^{*} S \Psi = I. 
\end{equation*}
We emphasize that the eigen decomposition is only used in the 
derivation of the method and is not performed in practical calculations.
For the moment we assume that $\lambda_1,\ldots,\lambda_N$ are positive and ordered from largest to smallest. From~\eqref{eqn:multishift} we have for each shift $z_l$
\begin{equation*}
	\begin{split}
	u_l &= (H - z_l S)^{-1} b \\
	&= \left(\Psi^{-*} \Lambda \Psi^{-1} - z_l \Psi^{-*}
	\Psi^{-1}\right)^{-1} b\\
	&= \Psi \left( \Lambda  - z_l \right)^{-1} \Psi^{*} b.
	\end{split}
\end{equation*}
Since  $\Re z\le 0,$ we may now use a pole expansion for $f^l(A) = 1/(A-z_lI)$ generated by the procedure in \cite{LinLuYingE2009} to write
\begin{equation}
	\begin{split}
		u_l^P & = \Psi \sum_{k=1}^{P}  \frac{\omega_{k}}{\xi_{k}-z_l}
		\left(\Lambda - \xi_{k}\right)^{-1} 
		\Psi^{*} b\\
		& = \sum_{k=1}^{P} \frac{\omega_{k}}{\xi_{k}-z_l}
		\left(\Psi^{-*}\Lambda\Psi^{-1} - \xi_{k}\Psi^{-*}\Psi^{-1}\right)^{-1} b \\
		&= \sum_{k=1}^{P} \frac{\omega_{k}}{\xi_{k}-z_l}
		\left(H - \xi_{k}S\right)^{-1} b, 
	\end{split}
	\label{eqn:pole_approx}
\end{equation}
which yields an approximate solution, $u_l^P,$ to the true solution $u_l.$ 
\begin{remark}
Using the pole expansion method to solve the parametrized systems~\eqref{eqn:multishift} motivated our choice of $f.$ However, the use of the pole
expansion only places mild requirements on $f$, so the method presented
here may potentially be used to solve systems with different types of
parametrization.
\end{remark}

To simplify the notation let us define $h_k$ as the solutions to the set of problems 
\begin{equation}
\label{eqn:pole_sub}
\left(H - \xi_{k}S\right)h_k  = b,
\end{equation} 
such that our approximate solution is now simply formed as
\begin{equation}
\label{eqn:solution}
u_l^P = \sum_{k=1}^{P} \frac{\omega_{k}}{\xi_{k}-z_l} h_k.
\end{equation}

If we assume that $h_k$ is computed exactly the solutions to~\eqref{eqn:pole_shift} formed via~\eqref{eqn:solution} satisfy the asymptotic, $P \rightarrow \infty,$ error bound:
\begin{equation*}
\frac{\|u_l - u_l^P\|_S}{\|\Psi^* b\|_2}  \leq \Or \left( e^{-C P / \log\left(\lambda_1/\lambda_N\right)} \right).
\end{equation*}
Given that the $h_k$ are computed exactly we may write
\begin{equation}
  \begin{split}
    \|u_l - u_l^P\|_S &= \|\Psi f^l(\Lambda) \Psi^* b - \Psi f^l_P(\Lambda) \Psi^* b\|_S \\
    &= \|f^l(\Lambda) \Psi^* b  - f^l_P(\Lambda) \Psi^* b\|_2 \\
    &\leq \|f^l(\Lambda)  - f^l_P(\Lambda) \|_2 \|\Psi^* b\|_2.
  \end{split}
\label{eqn:errS}
\end{equation}
Finally, using the error approximation~\eqref{eqn:HHTerror} in conjunction with~\eqref{eqn:errS} yields the desired result. Furthermore, in the special case where $S=I$ the error bound reduces to
\begin{equation*}
\frac{\|u_l - u_l^P\|_2}{\|b\|_2}  \leq \Or \left( e^{-C P / \log\left(\lambda_1/\lambda_N\right)} \right).
\end{equation*}

These bounds show that asymptotically the error decreases exponentially with respect to the number of poles used in the expansion for $f^l(A).$ However, there is additional error introduced since $h_k$ is computed inexactly. Thus, the overall error will often be dominated by the error in the computation of  $h_k.$ More specifically, let us define
\begin{equation}
\tilde{u}_l^P =  \sum_{k=1}^{P} \frac{\omega_{k}}{\xi_{k}-z_l} \tilde{h}_k
\label{eqn:approx_solution}
\end{equation}
where $\tilde{h}_k$ represents an approximation of $h_k.$ We assume that $\tilde{h}_k$ satisfies the error bound
\begin{equation}
\label{eqn:pole_sub_error}
\frac{\|\Psi^*(b - (H-\xi_kS)\tilde{h}_k)\|_2}{\|\Psi^*b\|_2} \leq \epsilon
\end{equation}
and define $r_k = b - (H-\xi_kS)\tilde{h}_k.$ In this case the solutions to~\eqref{eqn:pole_shift} formed via~\eqref{eqn:approx_solution} satisfy the asymptotic, $P \rightarrow \infty,$ error bound:
\begin{equation*}
\frac{\|u_l - \tilde{u}_l^P\|_S}{\|\Psi^*b\|_2} \leq  \Or \left( e^{-C P / \log\left(\lambda_1/\lambda_N\right)} \right) + \epsilon \left( \max_{i=1,\ldots,N}\frac{1}{\vert \lambda_i-z_l \vert} \right).
\end{equation*} 

Under the assumptions about the inexact solutions made here we may write
\begin{equation*}
  \begin{split}
    \tilde{u}_l^P &= \sum_{k=1}^{P} \frac{\omega_{k}}{\xi_{k}-z_l} \tilde{h}_k \\
                 &=  \sum_{k=1}^{P} \frac{\omega_{k}}{\xi_{k}-z_l} \left( h_k - (H-\xi_kS)^{-1}r_k \right) \\
                 &= u_l^P - \sum_{k=1}^{P} \frac{\omega_{k}}{\xi_{k}-z_l} (H-\xi_kS)^{-1}r_k.
  \end{split}
\end{equation*}
In conjunction with the error bound~\eqref{eqn:pole_sub_error} this implies that 
\begin{equation}
  \label{eqn:utilde_error}
  \begin{split}
    \frac{\|u_l - \tilde{u}_l^P\|_S}{\|\Psi^*b\|_2} &= \frac{\|u_l - u_l^P + \Psi f_P^l(\Lambda) \Psi^*r_k\|_S}{\|\Psi^*b\|_2} \\
             &\leq \frac{\|u_l - u_l^P\|_S}{\|\Psi^*b\|_2} + \frac{\|f_P^l(\Lambda) \Psi^*r_k\|_2}{\|\Psi^*b\|_2} \\
             &\leq \frac{\|u_l - u_l^P\|_S}{\|\Psi^*b\|_2} + \epsilon \|f_P^l(\Lambda)\|_2\\ 
             &\leq \frac{\|u_l - u_l^P\|_S}{\|\Psi^*b\|_2} + \epsilon \left( \|f_P^l(\Lambda) - f^l(\Lambda)\|_2 + \|f^l(\Lambda)\|_2\right)\\
             &\leq \frac{\|u_l - u_l^P\|_S}{\|\Psi^*b\|_2} + \epsilon\|f_P^l(\Lambda) - f^l(\Lambda)\|_2 + \epsilon \left( \max_{i=1,\ldots,N}\frac{1}{\vert \lambda_i-z_l \vert} \right).
  \end{split}
\end{equation}
Finally, using the estimate~\eqref{eqn:HHTerror} along with~\eqref{eqn:utilde_error} yields the desired result. Once again, in the case where $S=I$ the error bound simplifies to
\begin{equation*}
\frac{\|u_l - \tilde{u}_l^P\|_2}{\|b\|_2} \leq  \Or \left( e^{-C P / \log\left(\lambda_1/\lambda_N\right)} \right) + \epsilon \left( \max_{i=1,\ldots,N}\frac{1}{\vert \lambda_i-z_l\vert} \right).
\end{equation*} 

The error bound shows us that the error may be dominated by either the error in the pole expansion or the error in the solutions of~\eqref{eqn:pole_sub}. Since the error in the pole expansion decays exponentially, it is often best to control the overall error via the relative error requested when solving~\eqref{eqn:pole_sub}.  

Given that we are interested in solving systems of the form~\eqref{eqn:pole_shift} for a large number of shifts, the key observation
in~\eqref{eqn:pole_approx} is that the vectors $h_k$ are independent of
the shifts $z_l$ because the $\xi_k$ are independent of $z_l.$ Therefore, this method parametrizes the solutions to $N_z$ linear systems of
the form~\eqref{eqn:pole_shift} on the solutions of $P$  independent
linear systems of the form~\eqref{eqn:pole_shift}. In fact, once this parametrization has been done the method is completely flexible and any method may be used to solve the sub-problems. 

\begin{remark}
Because the $\xi_k$ appear in complex conjugate pairs, if $H,S$ and $b$ are real then so do the solutions to the systems~\eqref{eqn:pole_sub}. Therefore, in this situation we only have to solve $P/2$ systems, where for simplicity we assume $P$ is even.
\end{remark}

The bulk of the computational cost in this method is the
necessity of solving $P$ systems of the form~\eqref{eqn:pole_sub} after
which the vectors $h_k$ may be combined with different weights to yield
approximation solutions for as many distinct $z_l$ as desired. In fact, because the systems are completely independent this method
can be easily parallelized with up to $P$ (or $P/2$ depending on
whether symmetry is used) machines. Once the sub-problems have been solved, computing a solution for all $N_z$ shifts costs $O(PN_zN).$ 
Furthermore, as long as the $h_k$ are saved, solutions for new shifts
may be computed as needed with negligible computational cost. This is in
contrast to the Lanczos method where, as discussed in Section~\ref{sec:lanczos}, computing a solution for a new shift generally
requires running the algorithm again from the start unless all
$\{v_i\}$ are stored.

In some situations, the set of systems~\eqref{eqn:pole_sub} may even be
simultaneously solved using existing Krylov based methods for
simultaneously solving shifted systems, see, \eg,
\cite{Frommer2003,FrommerGlassner1998,SimonciniSzyld2007,DattaSaad1991,Meerbergen2003,GallivanGrimmeVan1996,BaiFreund2001,FeldmannFreund1995,SABAKI2012}.
Perhaps the simplest example would be to use the MINRES variation of
the method outlined in Section~\ref{sec:lanczos}. In this situation
the same methods may be applicable to the original systems, however,
the pole expansion reduces the number of shifts that have to be solved
for from $N_z$ to $P.$ Similarly, if a good preconditioner is known
for each of the $P$ distinct systems then each system may be
independently solved via an iterative method.

If the systems are amenable to the use of a direct method, \eg, LU
factorization or factorizations as in \cite{MARROK05}, then the $N_z$ shifted systems may be solved by computing
factorizations of the matrices $(H-\xi_kS)$ and using those
factorizations to solve the sub-problems. If a direct method is used the
procedure also allows for efficiently solving shifted systems with
multiple right hand sides. The specific solution methodologies we used
for the sub-problems of the pole expansion method will be discussed in
Section~\ref{sec:numer}.


\section{Indefinite Systems}\label{sec:indefinite}

Up until this point we have considered Hermitian matrix \response{pencils} $(H,S)$ for which the generalized eigenvalues are all positive. Motivated by the applications we discuss in Section~\ref{sec:connection} we now discuss the case where there are both positive and negative generalized eigenvalues and we seek a solution of a specific form. We still require that $H$ and $S$ are Hermitian and that $S$ is positive definite. We assume that the generalized eigenvalues of $(H,S)$ are ordered such that 
\begin{equation*}
  \lambda_1 \geq \dots \geq \lambda_M > 0 > \lambda_{M+1} \geq \dots \geq \lambda_N. 
\end{equation*}
Let $\Psi_{+} = [\psi_{1},\cdots,\psi_{M}]$ and $\Psi_{-} = [\psi_{M+1},\cdots,\psi_{N}].$ Similarly let $\Lambda_{+} =
\mathrm{diag}[\lambda_{1},\cdots,\lambda_{M}]$ and $\Lambda_{-} = \mathrm{diag}[\lambda_{M+1},\cdots,\lambda_{N}].$ 
While the notation used here mirrors that earlier in the paper, the generalized eigenvectors and eigenvalues here are distinct from the rest of the paper.

\subsection{Lanczos Method}
For the Lanczos based method used here there are two considerations that
have to be made with respect to indefinite systems. For the purposes of
this section we make the simplification, as in
Section~\ref{sec:lanczos}, that we first transform the problem in a
manner such that $S=I.$ Under this assumption we are interested in
solving systems for which the right hand side $\tilde{b}$ satisfies
$\Psi_{-}^*\tilde{b} = 0$ and the solution satisfies $\Psi_{-}^*x^l =
0.$ We note that in exact arithmetic $\Psi_{-}^*\tilde{b} = 0$ implies
that $\Psi_{-}^*x^l = 0$. However, we must ensure that the numerical
method used to solve these systems maintains this property.
Specifically, we need to ensure that $\Psi_{-}^*x^l_k \approx 0,$ and,
if we are using a CG style method, we must ensure that the sub-problems
remain non-singular. Both of these conditions are reliant upon the
mutual orthogonality between the Lanczos vectors and the columns of
$\Psi_{-}.$

In practice, we have observed that there is no excessive build up of
components in the $\Psi_{-}$ directions amongst the Lanczos vectors and
thus good orthogonality is maintained between the solutions we compute
using the Lanczos based methods and the negative generalized eigenspace.
Furthermore, in practice if we assume that $\Psi_{-}$ is known, then at
each step of the Lanczos process we may project out any components of
the Lanczos vectors that lie in the negative generalized eigenspace.
Such a procedure enforces orthogonality, up to the numerical error in
the projection operation, between the computed solution $x_k$ and
$\Psi_{-}.$

\subsection{\response{Pole expansion}}
The expansion we used in~\eqref{eqn:pole_approx} is valid for
Hermitian positive definite matrices. For
the case where the diagonal matrix of generalized eigenvalues $\Lambda$
has both positive and negative entries, 
the pole expansion does not directly give accurate results.
However, the pole expansion is
still applicable in the case where we wish to solve the systems
projected onto the positive generalized eigenspace. Section~\ref{sec:connection} provides the motivation for considering such
problems. Similar to before, we are interested in solving systems for
which the right hand side $b$ satisfies $\Psi_{-}^*b = 0$ and the
solution satisfies $\Psi_{-}^*Su_l = 0.$ To accomplish this we use a
pole expansion using $\lambda_1$ and $\lambda_M$ as the bounds of the
spectrum. We note that in exact arithmetic $\Psi_{-}^*b = 0$ implies
that $\Psi_{-}^*Su_l = 0.$ However we must ensure that our numerical
methods retain this property. 
Earlier in this section  we discussed
the impact of this generalization on the Lanczos based solver.
Here we restrict our discussion to the impact of solving an indefinite
system on the pole expansion method and provide an argument for why we do not observe difficulty in this regime.

If we once again assume that $\tilde{h}_k$ satisfies~\eqref{eqn:pole_sub_error} we may conclude that 
\begin{equation*}
  \frac{\|u_l - \tilde{u}_l^P\|_S}{\|\Psi_{+}^* b\|_2}  \leq \Or \left( e^{-C P / \log\left(\lambda_1/\lambda_M\right)} \right) + \Or \left( \epsilon \right) \left( \max_{i=1,\ldots,M}\frac{1}{\vert \lambda_i-z_l \vert} \right).
\end{equation*} 

To argue such a bound, in a manner similar to before, we may write
\begin{equation*}
  \begin{split}
    \frac{\|u_l - \tilde{u}_l^P\|_S}{\|\Psi_{+}^*b\|_2} &= \frac{\|u_l - u_l^P + \Psi f_P^l(\Lambda) \Psi^*r_k\|_S}{\|\Psi^*b\|_2} \\
             &\leq \frac{\|u_l - u_l^P\|_S}{\|\Psi^*b\|_2} + \frac{\|f_P^l(\Lambda) \Psi^*r_k\|_2}{\|\Psi_{+}^*b\|_2} \\
             &\leq \frac{\|u_l - u_l^P\|_S}{\|\Psi^*b\|_2} + \frac{\| f_P^l(\Lambda_{+})\Psi_{+}^*r_k\|_2 + \| f_P^l(\Lambda_{-})\Psi_{-}^*r_k\|_2 }{\|\Psi_{+}^*b\|_2} \\
             &\leq \frac{\|u_l - u_l^P\|_S}{\|\Psi^*b\|_2} + \epsilon \|f_P^l(\Lambda_{+})\|_2 + \frac{\| f_P^l(\Lambda_{-})\Psi_{-}^*r_k\|_2}{\|\Psi_{+}^*b\|_2}\\ 
             &\leq \frac{\|u_l - u_l^P\|_S}{\|\Psi^*b\|_2} + \epsilon \left( \|f_P^l(\Lambda_{+}) - f^l(\Lambda_{+})\|_2 + \|f^l(\Lambda_{+})\|_2\right) + \epsilon \| f_P^l(\Lambda_{-}) \|_2.
  \end{split}
\end{equation*}
Here we observe that the rational function $f_P^l$ is well behaved on the negative real axis in a manner dependent on the closest poles, which are of order $\lambda_M$ away. Furthermore, if the residuals $r_k$ are orthogonal to the negative general eigenspace this term vanishes. Also the construction of $f_P^l$ once again implies a dependence on the distance between the shifts and the generalized eigenvalues of the system. When using the pole expansion method on an indefinite system as long as the sub-problems are appropriately solved the use of the expansion basically maintains the accuracy of the overall solution. Furthermore, as long as the overall solution is computed accurately enough, and as long as $S$ is reasonably conditioned we cannot have large components of the computed solution in the negative generalized eigenspace. Furthermore, if necessary we may simply project out the components of $\tilde{h}_k$ in the negative generalized eigenspace to ensure that the overall solution error is not impacted by the sub-problem solution method.


\section{Connection with electronic structure calculation}\label{sec:connection}

In this section we discuss the connection between the problem with
multiple shifts in~\eqref{eqn:multishift} and several aspects
of the electronic structure theory, which are based on perturbative
treatment of Kohn-Sham density functional
theory~\cite{HohenbergKohn1964,KohnSham1965} (KSDFT). 
KSDFT is the most
widely used electronic structure theory for describing the ground state
electronic properties of molecules, solids and other nano structures. 
To simplify our discussion, we assume the computational domain is
$\Omega=[0,L]^3$ with periodic boundary conditions.  
We use linear algebra notation, and we do not distinguish integral
operators from their kernels.  For example, we may simply denote
$\hat{f}(\vr)=\hat{A}[\hat{g}](\vr)\equiv \int \hat{A}(\vr,\vr')
\hat{g}(\vr')\ud \vr'$ by $\hat{f}=\hat{A}\hat{g}$, and represent the
operator $\hat{A}$ by its kernel $\hat{A}(\vr,\vr')$. All
quantities represented in the real space are given in the form such as
$\hat{H}(\vr,\vr')$ and $\hat{f}(\vr)$, and the corresponding matrix or
vector coefficients represented in a finite dimensional basis set is
given in the form such as  $H$ and $f$.

The Kohn-Sham equation defines a nonlinear eigenvalue problem
\begin{equation}
  \begin{split}
    &\hat{H}[\hat{\rho}]\hat{\psi}_{i} = \varepsilon_{i} \hat{\psi}_{i},\\
    &\hat{\rho}(\vr) = \sum_{i=1}^{N_e} \abs{\hat{\psi}_{i}(\vr)}^2, \quad \int
    \hat{\psi}^{*}_{i}(\vr) \hat{\psi}_{j}(\vr) \ud \vr = \delta_{ij},
  \end{split}
  \label{eqn:KS}
\end{equation}
where $N_e$ is the number of electrons (spin degeneracy is omitted
here for simplicity). The eigenvalues
$\{\varepsilon_{i}\}$ are ordered non-decreasingly.  The lowest
$N_{e}$ eigenvalues $\{\varepsilon_{i}\}_{i=1}^{N_{e}}$ are called the occupied state energies, and 
$\{\varepsilon_{i}\}_{j>N_{e}}$ are called the unoccupied state
energies.  We assume
$\varepsilon_{N_{e}+1}-\varepsilon_{N_{e}}>0$, \ie\ the system is an
insulating system~\cite{Martin2004}.  
The eigenfunctions $\{\hat{\psi}_{i}\}_{i=1}^{N_e}$ define the
electron density $\hat{\rho}(\vr)$, which in turn defines the Kohn-Sham
Hamiltonian 
\begin{equation}
	\hat{H}[\hat{\rho}] = -\frac12 \Delta + \hat{V}_{c}[\hat{\rho}] + \Vxc[\hat{\rho}] + \Vion.
\label{eqn:ksham}
\end{equation}
Here $\Delta$ is the Laplacian operator for characterizing the kinetic
energy of electrons. 
\begin{equation*}
	\hat{V}_{c}[\hat{\rho}](\vr) \equiv
	\int \frac{\hat\rho(\vr')}{\abs{\vr-\vr'}} \ud \vr'
	\label{}
\end{equation*}
is the Coulomb potential which is linear with respect to the electron
density $\hat{\rho}$. 
$\Vxc[\hat{\rho}]$ is a nonlinear
functional of $\hat{\rho}$, characterizing the many body 
exchange and correlation effect. $\Vion$ is the electron-ion interaction
potential and is independent of $\hat{\rho}$.  Because the eigenvalue problem
(\ref{eqn:KS}) is nonlinear, it is often solved iteratively by a class
of algorithms called self-consistent field iterations
(SCF)~\cite{Martin2004}, until~\eqref{eqn:ksham} reaches self-consistency.

When the self-consistent solution of the Kohn-Sham equation is obtained,
one may perform post Kohn-Sham calculations for properties within and
beyond the ground state properties of the system.  Examples of such
calculations include the density functional perturbation theory
(DFPT)~\cite{BaroniGiannozziTesta1987,GonzeAllanTeter1992,BaroniGironcoliDalEtAl2001},
the GW
theory~\cite{Hedin1965,AryasetiawanGunnarsson1998,FriedrichSchindlmayr2006,UmariStenuitBaroni2010,PingRoccaGalli2013,GiustinoCohenLouie2010}
and the random phase approximation
(RPA) of the electron correlation energy~\cite{LangrethPerdew1975,LangrethPerdew1977,Furche2001,NguyenGironcoli2009}.
In these theories, a key quantity is the so called independent particle
polarizability matrix, often denoted by $\hat{\chi}_{0}(\omega)$.  The
independent particle polarizability matrix characterizes the first order
\textit{non-self-consistent} response of the electron density $\delta
\hat{\rho}_{0}(\omega)$ with respect to the time dependent external
perturbation potential $\delta \hat{V}(\omega)$, where $\omega$ is the
frequency of the time dependent perturbation potential.  $\omega$ can be
chosen to be $0$, characterizing the static linear response of the
electron density with respect to the static external perturbation
potential.  

In the density functional perturbation theory, the first order
\textit{self-consistent} static response (\ie\  the physical response) of
the electron density $\delta \hat{\rho}(\vr)$ with respect to the static
external perturbation potential $\delta \hat{V}$ can be computed as
\begin{equation*}
	\delta \hat{\rho} = \hat{\chi}(0) \delta \hat{V},
	\label{}
\end{equation*}
where the operator $\hat{\chi}(0)$ is directly related
to $\hat{\chi}_{0}(0)$ as 
\begin{equation*}
	\hat{\chi}(0) = \left[ I - \hat{\chi}_{0}(0) \left( \hat{V}_{c} + \frac{\delta
	\Vxc}{\delta \hat{\rho}} \right)
	\right]^{-1} \hat{\chi}_{0}(0).
	\label{}
\end{equation*}

Many body perturbation theories such as the GW
theory computes 
the quasi-particle energy which
characterizes the excited state energy spectrum of the system.
The key step for calculating the quasi-particle energy is the
computation of the screened Coulomb operator, which is defined as
\begin{equation*}
	\hat{W}(i\omega) = \left( I - \hat{V}_{c} \hat{\chi}_{0}(i\omega)
	\right)^{-1} \hat{V}_{c}.
	\label{}
\end{equation*}
Here $\chi_{0}(i\omega)$ should be computed on a large set of 
frequencies on the imaginary axis $i\omega$.

The random phase approximation (RPA) of the electron correlation energy
improves the accuracy of many existing exchange-correlation functional in
ground state electronic structure calculation.
Using the adiabatic connection
formula~\cite{LangrethPerdew1975,LangrethPerdew1977}, the correlation
energy can be expressed as
\begin{equation*}
	E_{c} = -\frac{1}{2\pi} \int_{0}^{1} \int_{0}^{\infty} 
	\Tr\left\{ \hat{V}_{c} [\hat{\chi}_{\lambda}(i\omega) -
	\hat{\chi}_{0}(i\omega)] \right\}
	\ud\omega \ud \lambda,
	\label{}
\end{equation*}
and $\hat{\chi}_{\lambda}(i\omega)$ can be computed from
$\hat{\chi}_{0}(i\omega)$ as
\begin{equation*}
	\hat{\chi}_{\lambda}(i\omega)	 = \hat{\chi}_{0}(i\omega) + \lambda
	\hat{\chi}_{0}(i\omega)
	\hat{V}_{c} \hat{\chi}_{\lambda}(i\omega).
	\label{}
\end{equation*}

In all the examples above, computing $\hat{\chi}_{0}(i\omega)$ is
usually the bottleneck.  For simplicity we consider the case 
when $\hat{H}$ is real and therefore the eigenfunctions
$\hat{\psi}_{i}(\vr)$ are also real. In such case the kernel
$\hat{\chi}_{0}(i\omega)(\vr,\vr')$ can be computed from the Adler-Wiser
formula~\cite{Adler1962,Wiser1963} as
\begin{equation*}
	\hat{\chi}_{0}(\vr,\vr',i\omega) = 2 \Re \sum_{i \le
	N_{e}, j > N_{e}} 
	\frac{\hat{\psi}_{i}(\vr)\hat{\psi}_{j}(\vr)\hat{\psi}_{j}(\vr')\hat{\psi}_{i}(\vr')}
	{\varepsilon_{i}-\varepsilon_{j}+i\omega} .
\end{equation*}
The summation $\sum_{j > N_{e}}$ requires the computation of a
large number of eigenstates which is usually prohibitively expansive.  
Recent
techniques~\cite{UmariStenuitBaroni2010,PingRoccaGalli2013,GiustinoCohenLouie2010}
have allowed to avoid the direct computation of
$\{\hat{\psi}_{j}\}_{j>N_{e}}$ when $\hat{\chi}_{0}(i\omega)$  is
multiplied to an arbitrary vector $\hat{g}$ as
\begin{equation}
	\begin{split}
		\hat{\chi}_{0}(i\omega)[\hat{g}](\vr) &=	2 \Re \sum_{i \le N_{e}} 
		\hat{\psi}_{i}(\vr) \int 
		\sum_{j > N_{e}}
		\frac{\hat{\psi}_{j}(\vr)\hat{\psi}_{j}(\vr')}{\varepsilon_{i}-\varepsilon_{j}
		+ i\omega } 
		\hat{\psi}_{i}(\vr') \hat{g}(\vr')\ud \vr'\\
		& \equiv 2 \Re \sum_{i \le N_{e}} 
		\hat{\psi}_{i}(\vr) \hat{u}_{i}(\vr).
	\end{split}
	\label{eqn:chi0}
\end{equation}
Here $\hat{u}_{i}(\vr)$ can be solved through the equation
\begin{equation}
	\hat{Q} \left[\hat{H} - (\varepsilon_{i}+i\omega)\right] \hat{Q}
	\hat{u}_{i} =
	-\hat{Q}[\hat{\psi}_{i}\odot \hat{g}].
	\label{eqn:multishiftOriginal}
\end{equation}
The operator $\hat{Q}(\vr,\vr')=\delta(\vr,\vr') - \sum_{i\le
N_{e}}\hat{\psi}_{i}(\vr)\hat{\psi}^{*}_{i}(\vr')$ is a projection
operator onto the occupied states (noting
$\hat{\psi}^{*}_{i}(\vr')=\hat{\psi}_{i}(\vr')$ is real), and
$[\hat{\psi}_{i} \odot \hat{g}](\vr) \equiv
\hat{\psi}_{i}(\vr)\hat{g}(\vr)$ is the element wise product between two
vectors.

Without loss of generality we may set the largest occupied state energy $\varepsilon_{N_{e}}=0$.
Equation~\eqref{eqn:multishiftOriginal} can be
reduced to the form
\begin{equation}
	(\hat{H} - z) \hat{u} = \hat{f},
	\label{eqn:multishiftFine}
\end{equation}
for multiple shifts $z = \varepsilon_{i} + i\omega$ ($\Re z \le 0$) and multiple right
hand sides $\hat{f}$.  Furthermore, $\hat{u},\hat{f} \in \Ran(
\hat{Q})$ where $\Ran$ is the range of the operator $\hat{Q}$.  
Equation~\eqref{eqn:multishiftFine} can be solved in practice using a finite
dimensional basis set, such as finite element, plane waves, or
more complicated basis functions such as numerical atomic
orbitals~\cite{SolerArtachoGaleEtAl2002}.
We denote the basis set by a collection of column vectors as 
$\hat{\Phi}=[\hat{\varphi}_{1}(\vr),\cdots,\hat{\varphi}_{N}(\vr)]$.
The overlap matrix associated with the basis set $\hat{\Phi}$ is
\begin{equation*}
	S = \hat{\Phi}^{*} \hat{\Phi},
	\label{}
\end{equation*}
and the projected Hamiltonian matrix in the basis $\hat{\Phi}$ is 
\begin{equation*}
	H = \hat{\Phi}^{*} \hat{H} \hat{\Phi}.
	\label{}
\end{equation*}
Using the ansatz that both the solution and the right hand side can be
represented using the basis set $\hat{\Phi}$ as
\begin{equation}
	\hat{u} = \hat{\Phi} u,\quad \hat{f} = \hat{\Phi} f,
	\label{eqn:ufbasis}
\end{equation}
\eqref{eqn:multishiftFine} becomes
\begin{equation*}
	(H - z S) u = \hat{\Phi}^* \hat{f} = S f,
	\label{}
\end{equation*}
which is~\eqref{eqn:multishift} with $b=S f$. Using the
eigen decomposition of the matrix \response{pencil} $(H,S)$ as in~\eqref{eqn:GE}, each eigenfunction in the real space
$\hat{\psi}_{i}(\vr)$ is given using the basis set
$\hat{\Phi}$ as 
\begin{equation}
	\hat{\psi}_{i} = \hat{\Phi}\psi_{i}.
	\label{eqn:psibasis}
\end{equation}
Combining Eqs.~\eqref{eqn:psibasis} and~\eqref{eqn:ufbasis}, the
condition $\hat{u}\in \Ran(\hat{Q})$ becomes
\begin{equation}
	\hat{\psi}_{i}^{*} \hat{u} = \psi_{i}^{*} S u = 0,\quad \forall i\le
	N_{e},
	\label{eqn:ucond}
\end{equation}
and similarly $\hat{f}\in \Ran(\hat{Q})$ becomes
\begin{equation}
	\psi_{i}^{*} S f = 0,\quad \forall i\le N_{e},
	\label{eqn:fcond}
\end{equation}

In practice the conditions~\eqref{eqn:ucond} and~\eqref{eqn:fcond} can
be satisfied by a projection procedure as described in
Section~\ref{sec:indefinite}.


\section{Numerical results}\label{sec:numer}
First we illustrate the accuracy and the scaling of the pole expansion
method.
We then present two distinct
numerical examples based on the general method presented here. In one
example we consider the case where $S=I$ and an iterative method is used
to solve the sub-problems, and in the second example we consider the
case where $S \neq I$ and a direct method is used to solve the
sub-problems. In each case the method presented here is compared with
the Lanczos style method described in Section~\ref{sec:lanczos}.

All of the numerical experiments were run in MATLAB on a Linux machine
with four 2.0GHz eight core CPUs and 256GB of RAM.  


\subsection{Accuracy of the pole expansion} 
\begin{figure}
  \centering
  \begin{subfigure}[b]{.49\textwidth}
    \includegraphics[width = 1\textwidth]{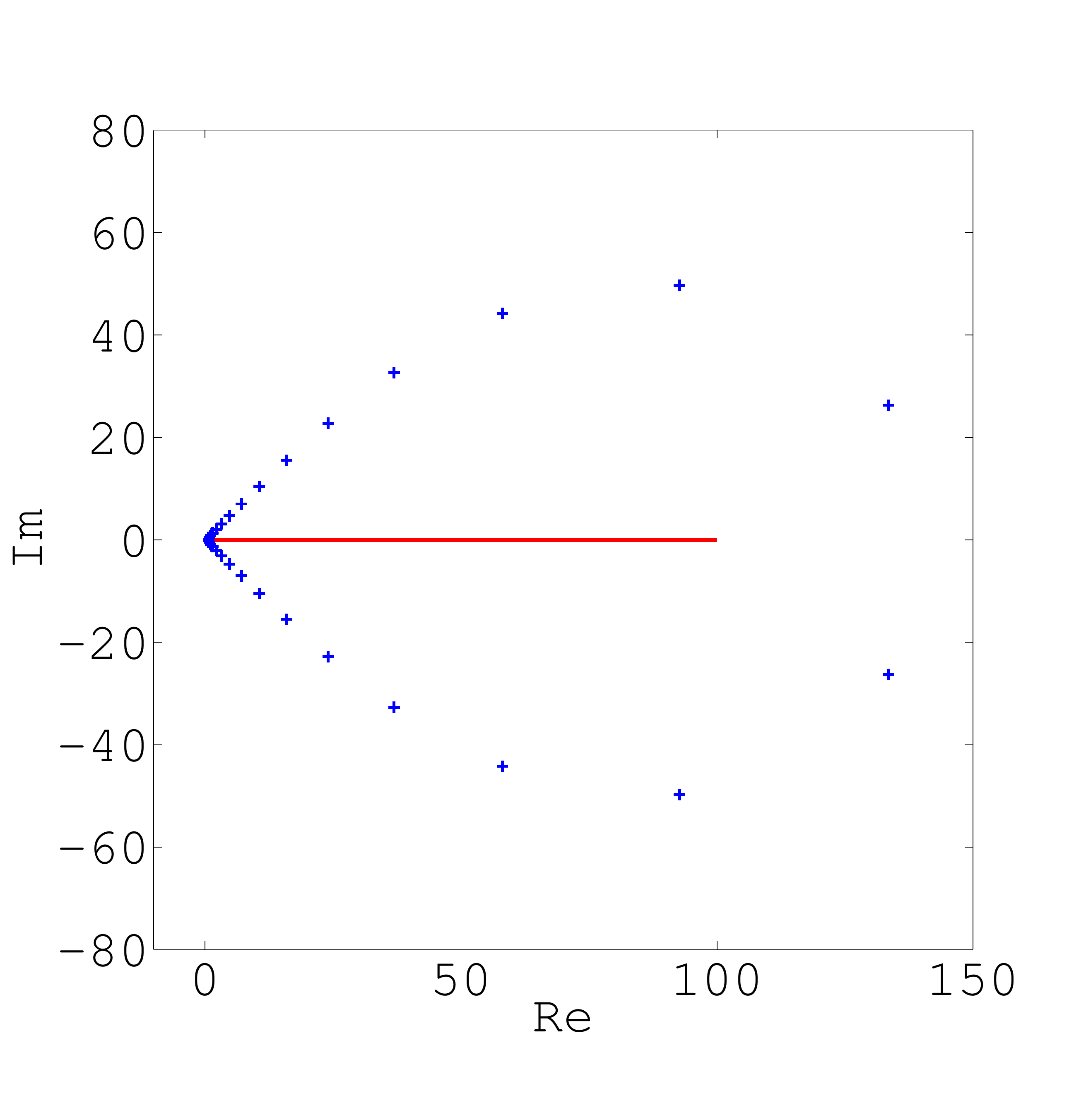} 
    \caption{View of the entire contour.}
  \end{subfigure}
  \begin{subfigure}[b]{.49\textwidth}
    \includegraphics[width = 1\textwidth]{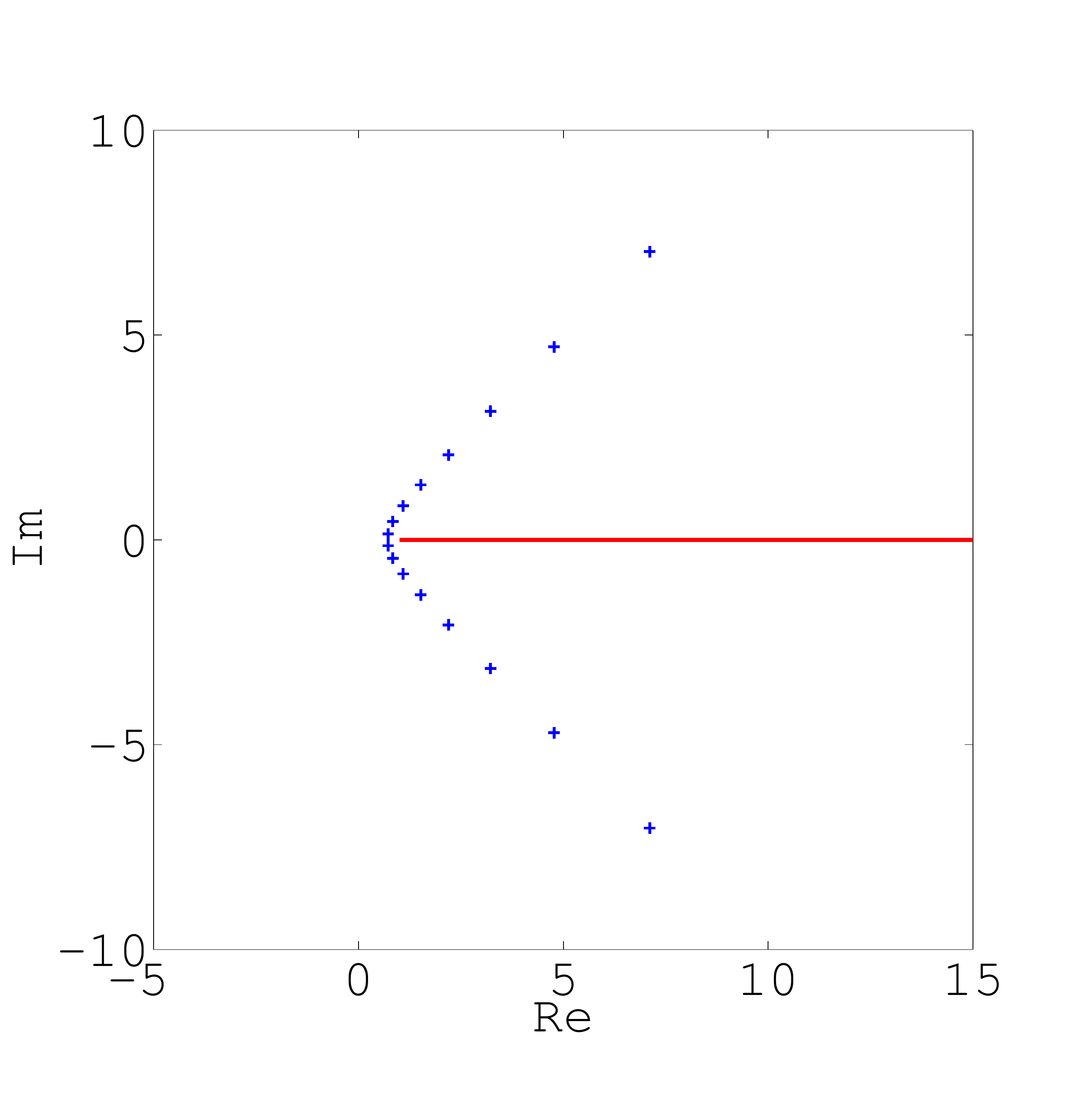} 
    \caption{View near the origin.}
  \end{subfigure}
  \caption{Example of quadrature nodes (blue $+$'s) used when the spectrum of $A$ (red line) lies in the interval $[1,100].$}
    \label{fig:contour_ex} 
\end{figure}

Before presenting the examples motivated by the preceding section, we
first demonstrate the behavior of the pole expansion method for 
approximating 
\begin{equation*}
f(x;z) = \frac{1}{x-z}
\end{equation*}
on some interval of the positive real axis.
Fig.~\ref{fig:contour_ex}
shows an example contour as computed via the method in Section~\ref{sec:pole} when the region of interest is $[1,100]$.

\begin{figure}
  \centering
  \begin{subfigure}[b]{.49\textwidth}
    \includegraphics[width = 1\textwidth]{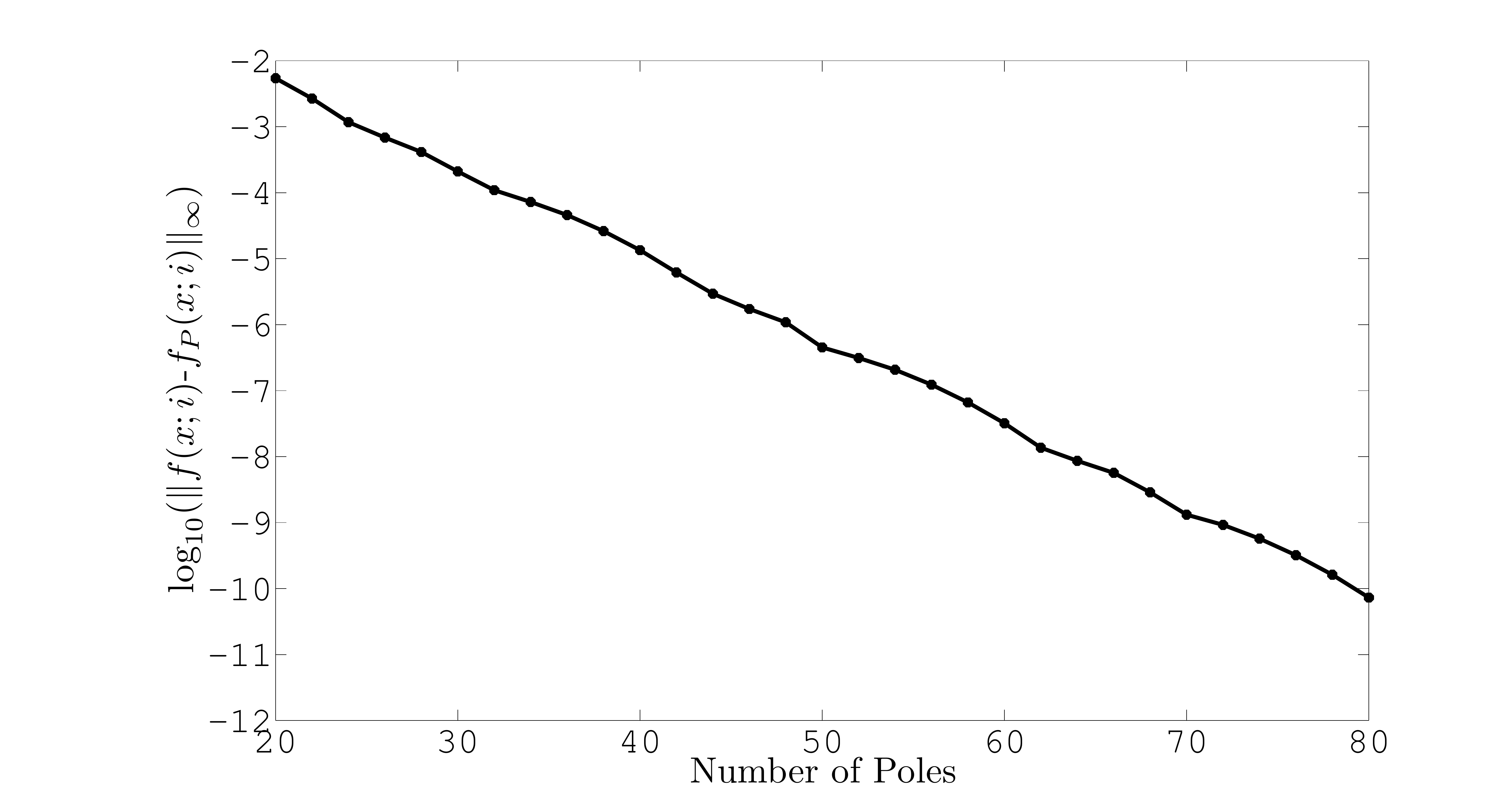} 
    \caption{}
    \label{fig:pole_err}  
  \end{subfigure}
  \begin{subfigure}[b]{.49\textwidth}
    \includegraphics[width = 1\textwidth]{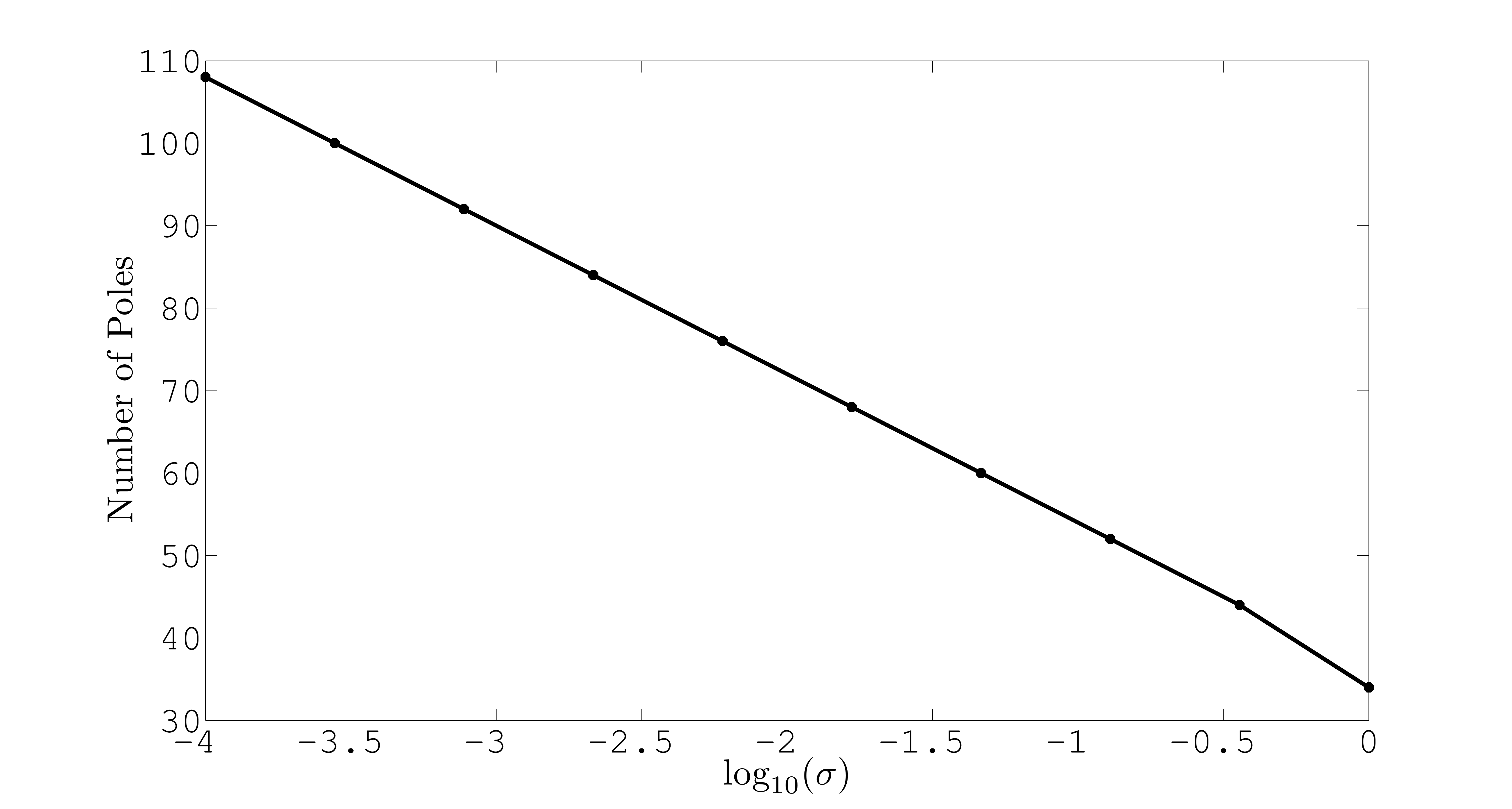} 
    \caption{}
    \label{fig:numpoles}
  \end{subfigure}
  \caption{Behavior of the pole expansion. (a) $\|f_P(x;i)-f(x;i)\|_{\infty},$ on a $\log_{10}$ scale, for the interval $x\in [1,1000]$ as the number of poles, $P,$ is increased. (b) Number of poles required to achieve $\|f_P(x;i)-f(x;i)\|_{\infty} \leq 10^{-8}$ for $x\in [\sigma,10].$}
\end{figure}

First we consider approximating $f(x;z)$ for a fixed value $z=i$ on the
interval $[1,1000].$
Fig.~\ref{fig:pole_err} shows the
$\|\cdot \|_{\infty}$ error of approximating $f(x;i)$ by $f_P(x;i),$
computed via sampling at $10000$ equally spaced points in $[1,1000]$
along the $x$ direction, as
the number of poles $P$ increases. We observe the 
exponential decay in error as the number of poles increases, which is aligned with the error analysis presented earlier. 

Next we illustrate the number of poles required to reach fixed accuracy 
as $x$ approaches $0$, \ie\ as $E_{g}\to 0$. 
We consider approximating $f(x;i)$ via $f_P(x;i)$ on the interval $[\sigma,10]$
for 20 
values of $\sigma \in [10^{-4},1]$ that are equally spaced on the
logarithmic scale.
Fig.~\ref{fig:numpoles} shows the number of poles required such that
$\|f_P(x;i)-f(x;i)\|_{\infty}$ is less than $10^{-8}$ over the interval
$[\sigma,10].$

\begin{figure}
  \centering
  \includegraphics[width = .75\textwidth]{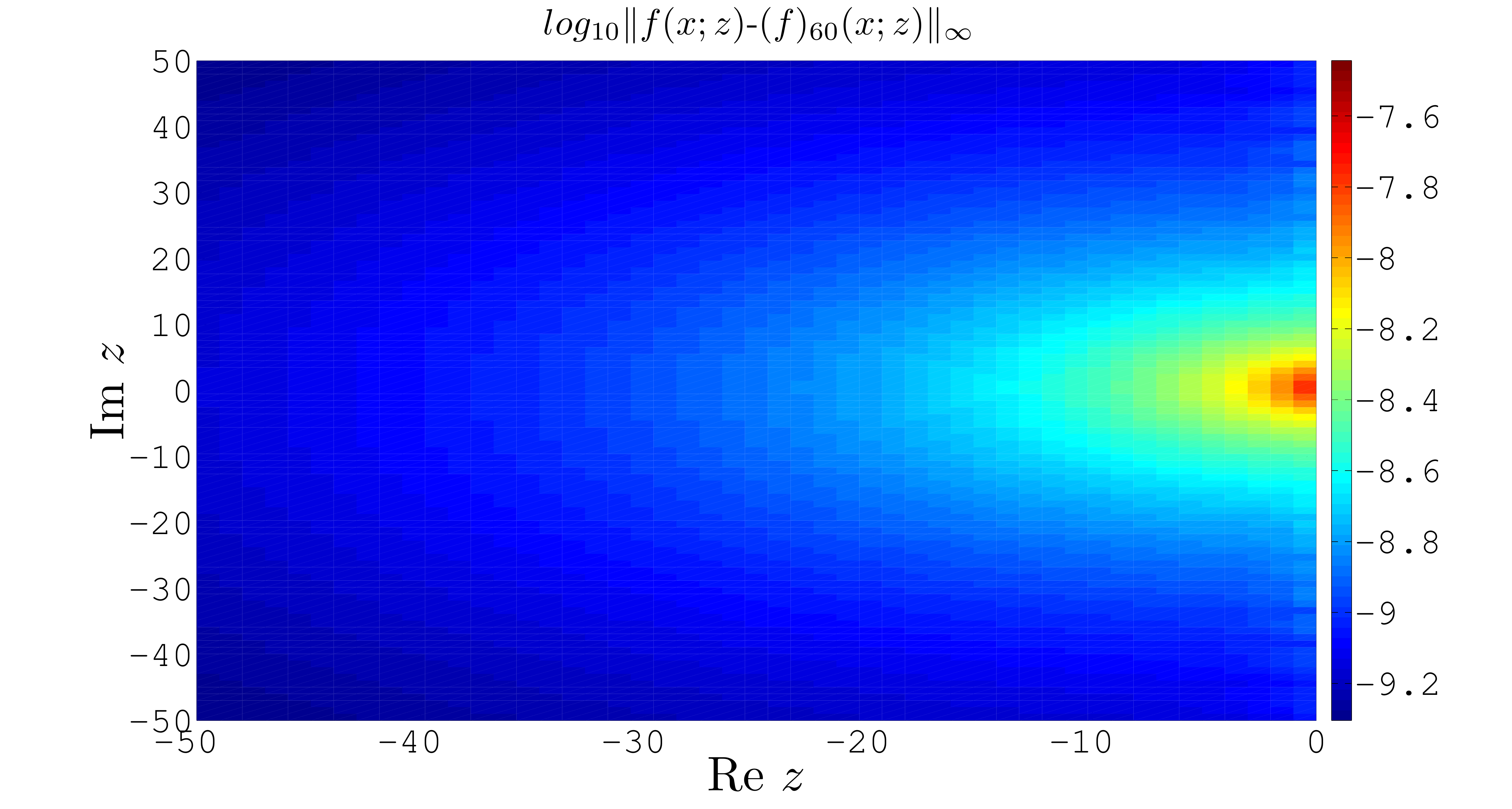} 
  \caption{$\|f_P(x;z)-f(x;z)\|_{\infty}$ on a $\log_{10}$ scale for the interval $x\in [1,1000]$ as $z$ is varied.}
  \label{fig:pole_range}  
\end{figure}

Finally we consider the approximation of $f(x;z)$ for a wide range
of $z$ with $\Re z \leq 0$, and demonstrate how $\| f_P(x;z) - f(x;z)
\|_{\infty}$ varies as $z$ changes  using the same pole expansion. We
fixed the number of poles used in the approximation to be 60. Fig.~
\ref{fig:pole_range} shows that high accuracy is maintained for all
$z$ in the left half plane. Here we kept the region of interest as $x\in [1,1000]$ and
consider the accuracy at 5000 distinct $z$ evenly distributed with $\Re z \in [-50,0]$ and $\Im z \in
[-50,50].$  

\subsection{Orthogonal basis functions: $S=I$}
First, we consider the case where orthogonal basis functions are used
and thus in the notation here $S=I.$ The Hamiltonian matrix $H$ takes
the form
\begin{equation*}
H = -\frac{1}{2}\Delta + V,
\end{equation*}
and the operator $V$ is obtained by solving the Kohn-Sham density
functional theory problem for a benzene molecule using the KSSOLV
package~\cite{YangMezaLeeEtAl2009}, which is a MATLAB toolbox for
solving Kohn-Sham equations for small molecules and solids implemented
entirely in MATLAB m-files.  The benzene molecule has $30$ electrons and
$15$ occupied states (spin degeneracy of $2$ is counted here). 
The atomic configuration of the benzene molecule is given in
Fig.~\ref{fig:benzene}.

The computational domain is $[0,20) \times [0,20)
\times [0,10)$ along the $x,y,z$ directions, respectively.  The
computational domain is discretized into $64\times 64 \times 32 =
131072$ points.  The Laplacian operator is discretized using the plane
wave basis set, and this set of grid corresponds to the usage of the kinetic
energy cutoff at $50.5$ Hartree. The Laplacian operator is applied using
spectral method, and is done efficiently using the Fast Fourier
Transform (FFT). The $15$ negative eigenvalues (occupied states) and the
corresponding eigenvectors are 
computed using 
the locally optimal block preconditioned conjugate gradient (LOBPCG)
method~\cite{Knyazev2001},
with a preconditioner of the form $(-\frac{1}{2}\Delta+.001)^{-1}.$
We only use these computed eigenvectors to ensure that the
right hand side for the set of equations we solve is orthogonal to the
negative eigenspace of the $(H,S)$ \response{pencil}, as required in
Section~\ref{sec:indefinite} for using the pole expansion for indefinite
systems.


We now solve the set of problems 
\begin{equation}
\label{eqn:n1}
(H-i\eta_lI)x_l = b,
\end{equation}
for $101$ equispaced $\eta_l$ in the interval $[-10,10].$ 
Because $H$ and $b$ are real, some systems are in essence solved
redundantly. However, here we are interested in the performance given the
number of shifts and not the specific solutions.  


We monitor the behavior of our method and the Lanczos method as the
condition number of $H$ increases. To this end we sequentially refine
the number of discretization points in each direction by a factor of $2$
and generate a potential function $V$  via Fourier interpolation. The
grid size is denoted by $N\times N\times N/2$, and the largest problem
we consider is discretized on a $256\times 256 \times 128$ grid. For
these large problems the eigenvectors associated with negative
eigenvalues are approximated by Fourier interpolates of the computed
eigenvectors for the smallest problem. 


To compare the methods we solved the set of
problems~\eqref{eqn:n1} via the pole expansion method using a number of
poles that depended on the size of the problem. To combat the slight
loss of accuracy that occurs for a fixed number of poles as the
condition number of the matrix grows, we increased the number of poles as
the problem size grew. During this step preconditioned
GMRES \cite{Saad1986} in MATLAB
was used to solve the sub-problems associated with the pole expansion.
The preconditioner used \response{is} of the form
\begin{equation*}
\left( -\frac{1}{2}\Delta + \xi_k\right)^{-1},
\end{equation*}
and similarly to $H$ it was efficiently applied via the FFT. The
requested accuracy of the GMRES routine is that the relative residual
is less than $10^{-7}.$ Since $H$ and $b$ are real only $P/2$
sub-problems had
to be solved. We remark that though the solution for different poles can
be straightforwardly parallelized, here we performed the calculations sequentially 
in order to compare with the sequential
implementation of the Lanczos method. Let $\tilde{x}_l^P$ denote the
approximate solutions computed using the pole expansion. The relative
error metric
\begin{equation*}
r_l^P = \frac{ \| b - (H - i\eta_l I )\tilde{x}_l^P \|_2 }{\| b \|_2}
\end{equation*}
is then computed for each shift after the approximate solutions have been computed. Finally, the Lanczos procedure described in Section~\ref{sec:lanczos} is called with a requested error tolerance of 
\begin{equation}
\label{eqn:lanczos_stop}
\frac{ \| b - (H - i\eta_l I )\hat{x}_l\|_2 }{\| b \|_2} \leq \max_{l} r_l^P,
\end{equation}
where $\hat{x}_l$ denotes the approximate solution computed via the Lanczos method. Since the Lanczos method reveals the residual for each shift at each iteration this stopping criteria is cheap computationally. Thus, the Lanczos method was run until all of the approximate solutions met the stopping criteria~\eqref{eqn:lanczos_stop}. The Lanczos method used here takes advantage of the cheap updates briefly described in Section~\ref{sec:lanczos} and to further save on computational time, once a solution for a given $\eta_l$ was accurate enough the implementation stopped updating that solution. For comparison purposes we also solved the problem with a version of the Lanczos algorithm that uses Householder reflectors to maintain orthogonality amongst the Lanczos vectors \cite{GVL}. For the smallest size problem used here, and without any shifts, the version from Section~\ref{sec:lanczos} took 460 iterations to converge to $10^{-8}$ accuracy while the version that maintained full orthogonality took 456 iterations to converge to the same accuracy. However, the method that maintained full orthogonality took 65 times longer to run and would be prohibitively expensive for the larger problems given the increased problem size and iteration count so for all the comparisons here we use the CG style method outlined in Section~\ref{sec:lanczos}.

The spectrum of the operator $H$ grows as $\mathcal{O}(N^2)$, and we
observe that the number of iterations 
required for the Lanczos method to converge grow roughly as $\mathcal{O}(N)$. We
do not expect to see such growth of the number of
iterations in the pole expansion method, given the preconditioner used
in solving the sub-problems. Furthermore, even though the cost per
iteration of solving for multiple shifts scales linearly in the Lanczos
method, the time required is dependent on the number of iterations
required to converge.  In fact, in this case where $H$ may be applied
very efficiently for a very large number of shifts, the cost at each
iteration may be dominated by the additional computational cost
associated with each shift. In contrast, once the sub-problems have been
solved, the pole expansion method has a fixed cost for computing the
solutions for all the $\eta_l$, which only depends on the number of
shifts and the number of poles.  If $\vert \eta_l \vert$ is large the
Lanczos method will converge very quickly for this specific problem, so
in the case where only shifts with large magnitude are considered the
Lanczos method may perform better. However, Section~\ref{sec:connection}
motivates our use of shifts spaced out along a portion of the imaginary
axis that includes 0.

Table~\ref{tab:numer1} reports the results of the pole
expansion method and the Lanczos method for the problem
described above. For the pole expansion method, the number of iterations
reported is the total number of iterations required to solve all of the
sub-problems. In both cases the error reported is the maximum computed
relative residual over the shifts. Finally, the total time taken to
solve the problems is reported for each method. 

We observe that as expected the Lanczos method performs better for
the smallest sized problem. However once we reach the mid sized problem
the methods perform comparably, and the pole expansion method is
actually a bit faster even though it takes a few more iterations
overall. What is important to notice is that the overall iteration count
of solving all the sub-problems associated with the pole expansion
method remains relatively constant even with the increased number of
poles. By the time we reach the largest
problem the pole expansion method outperforms the Lanczos
method, taking about half the time to solve the set of problems.

\begin{remark}
Because the GMRES method used here uses Householder transforms to
maintain orthogonality amongst the Krylov basis it is not very memory
efficient. However, we also ran this example using preconditioned TFQMR \cite{TFQMR}
in place of GMRES, and while the iteration count more than doubled for
the pole expansion method it was actually faster since the applications
of the Householder matrices in GMRES is very expensive. This comparison
actually helps demonstrate the flexibility of the method with respect to
the solver used. If a fast method was not available for applying $H$
then GMRES may be preferable, however in the case where $H$ may be
applied very efficiently such as TFQMR may be better suited to
the problem.
\end{remark}

\begin{table}
  \centering
  \begin{tabular}{c | c | c | c | c | c | c | c | }
    \cline{2-8}
    & \multicolumn{3}{|c|}{Lanczos} & \multicolumn{4}{|c|}{Pole Expansion} \\ \hline
    \multicolumn{1}{|c|}{Problem Size} & Iter. & Time(s) & Error & P & Iter. & Time(s) & Error \\ \hline
    \multicolumn{1}{|c|}{$64\times 64 \times 32$} & 297  & 97.23  & 9.83 $\times 10^{-7}$ & 70 & 1005 & 181.99 & 9.85 $\times 10^{-7}$  \\ \hline
    \multicolumn{1}{|c|}{$128\times 128 \times 64$} & 911  & 2073.89  & 6.91 $\times 10^{-6}$  & 80 & 969 & 1920.75  & 6.93 $\times 10^{-6}$ \\ \hline
    \multicolumn{1}{|c|}{$256\times 256 \times 128$} & 1746  & 35129.79  & 6.83 $\times 10^{-6}$  & 90 & 959 & 17127.93  & 6.85 $\times 10^{-6}$ \\ \hline
  \end{tabular}
  \caption{Comparison of the pole expansion method and the Lanczos method for solving $(H-i\eta_lI)x_l = b.$}
  \label{tab:numer1}
\end{table}

\begin{figure}
  \centering
  \begin{subfigure}[b]{.49\textwidth}
    \includegraphics[width = 1\textwidth]{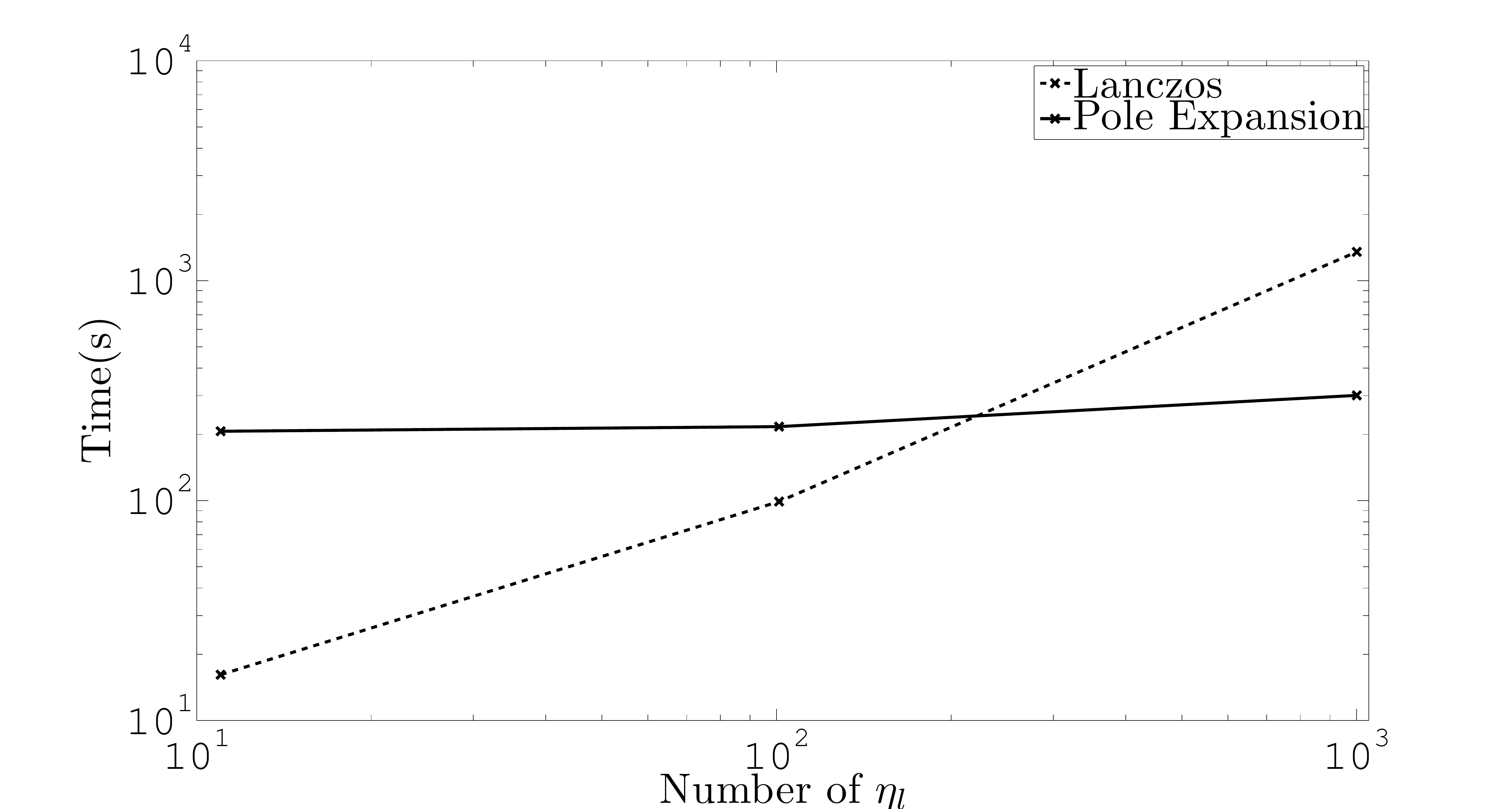} 
    \caption{$64\times 64 \times 32$}
    \label{fig:omega_small}  
  \end{subfigure}
  \begin{subfigure}[b]{.49\textwidth}
    \includegraphics[width = 1\textwidth]{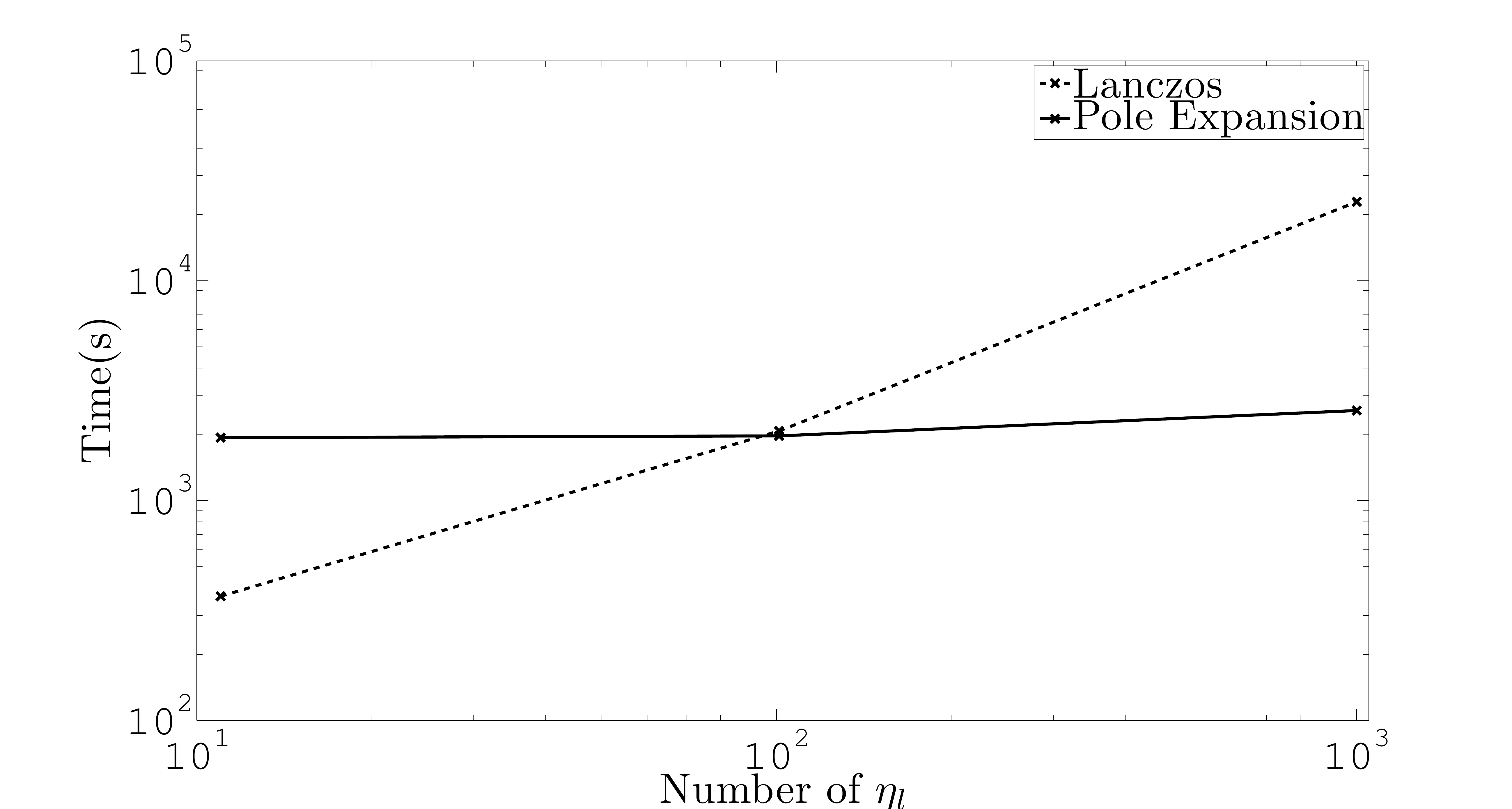} 
    \caption{$128\times 128 \times 64$}
    \label{fig:omega_large}  
  \end{subfigure}
   \caption{Time taken to solve problems for a varying number of $\eta_l$ equispaced in $[-10,10]$}
\end{figure}

To further demonstrate the scaling of the methods with respect to the
number of shifts we ran the problem at a fixed size and varied the
number of shifts. As before, preconditioned GMRES was used in the pole
expansion method and no parallelism was used when solving the
sub-problems. The same strategy as above was used to ensure that the
Lanczos method stopped once it had solved all the problems as accurately
as the least accurate solution found using the pole expansion. Fig.~
\ref{fig:omega_small} shows the time taken to solve the problems of size
$64\times 64 \times 32$ for a varying number of $\eta_l$ equispaced in
$[-10,10].$ In all cases the largest relative residual was on the order
of $7\times 10^{-7}.$ Similarly, Fig.~\ref{fig:omega_large} shows the
time taken to solve the problems of size $128\times 128 \times 64$ for a
varying number of $\eta_l$ equispaced in $[-10,10].$ In all cases the
largest relative residual was on the order of $7\times 10^{-6}.$ Here we
observe that in both cases the pole expansion method scales very well as
the number of shifts increases, especially in comparison to the Lanczos
method. For example, even in the case where the problem is $64\times
64\times 32$ and the Lanczos method takes around a third of the
iterations of the pole expansion method, if you take the number of
shifts to be large enough the pole expansion method becomes considerably
faster.

\subsection{Computation of $\hat{\chi}_{0}(i\omega)[\hat{g}](\vr)$}
Motivated by our discussion in Section~\ref{sec:connection} we used our
technique to compute $\hat{\chi}_{0}(i\omega)[\hat{g}](\vr)$ as defined
in~\eqref{eqn:chi0} where we only need $\hat{\psi}_i(\vr)$ for $i\leq
N_e$. In order to compute $\hat{\chi}_{0}(i\omega)[\hat{g}](\vr)$ for a
large number of $i\omega$ we had to compute $\hat{u}_i(\vr)$ via
\eqref{eqn:multishiftOriginal}.
Physically this corresponds to the calculation of the response of the
electron density with respect to external perturbation potential
$\hat{g}(\vr)$.

\begin{figure}
  \centering
  \begin{subfigure}[b]{.3\textwidth}
    \includegraphics[width = 1\textwidth]{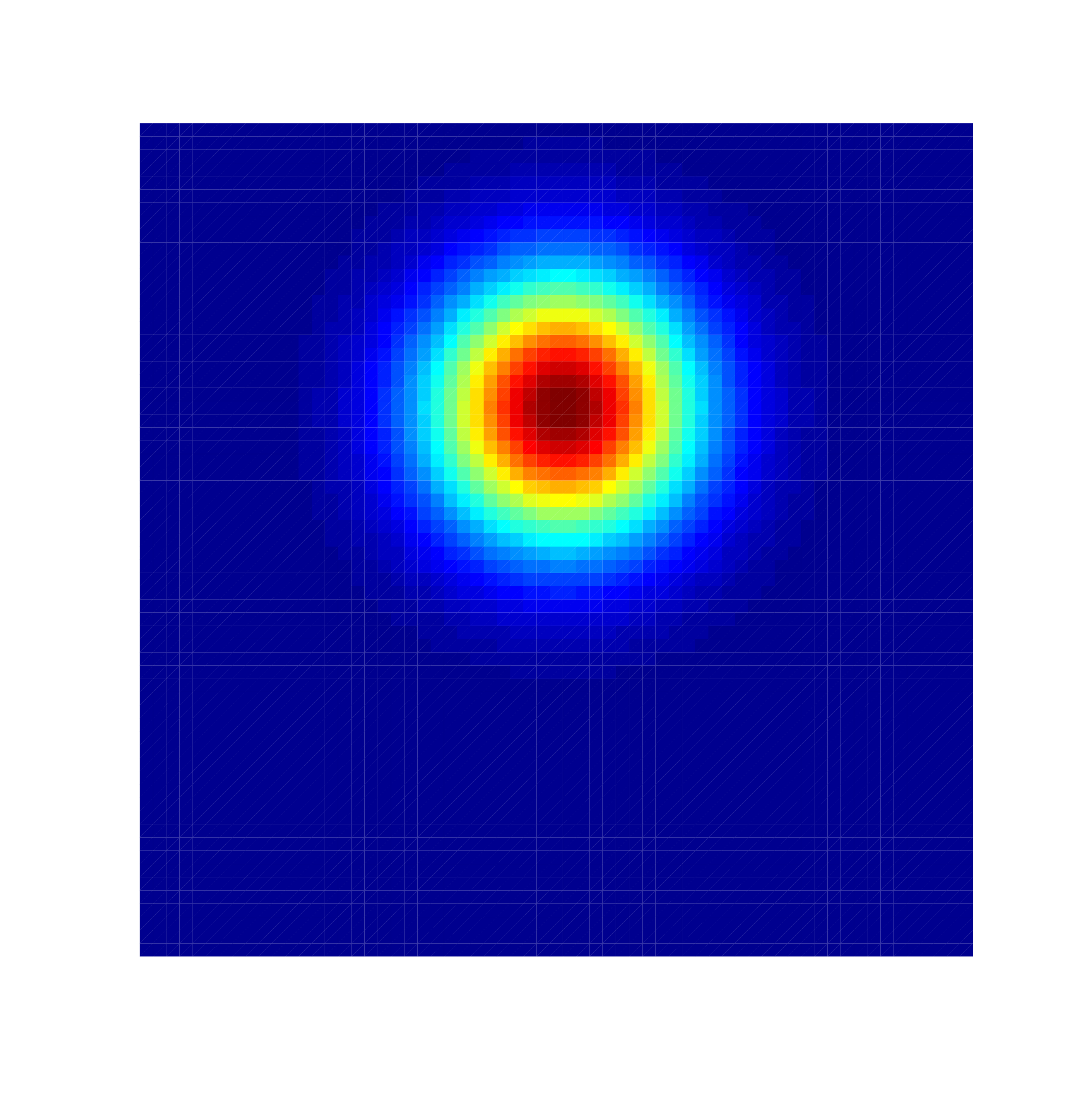} 
		\caption{Slice of  $\hat{g}.$}
    \label{fig:g}  
  \end{subfigure}
   \begin{subfigure}[b]{.3\textwidth}
    \includegraphics[width = 1\textwidth]{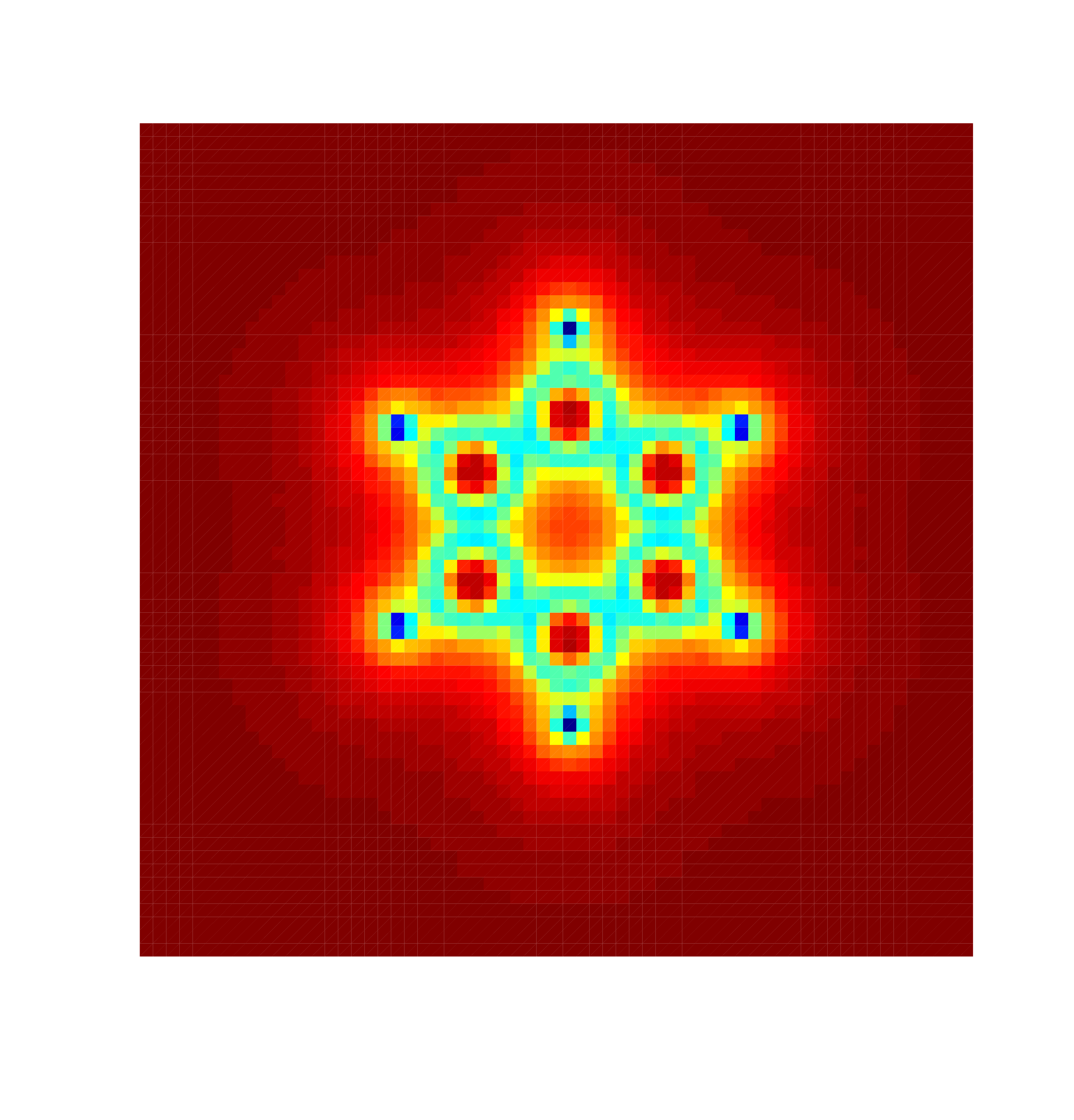} 
    \caption{Slice of  $V.$}
    \label{fig:V}  
  \end{subfigure}
  \begin{subfigure}[b]{.3\textwidth}
    \includegraphics[width = 1\textwidth]{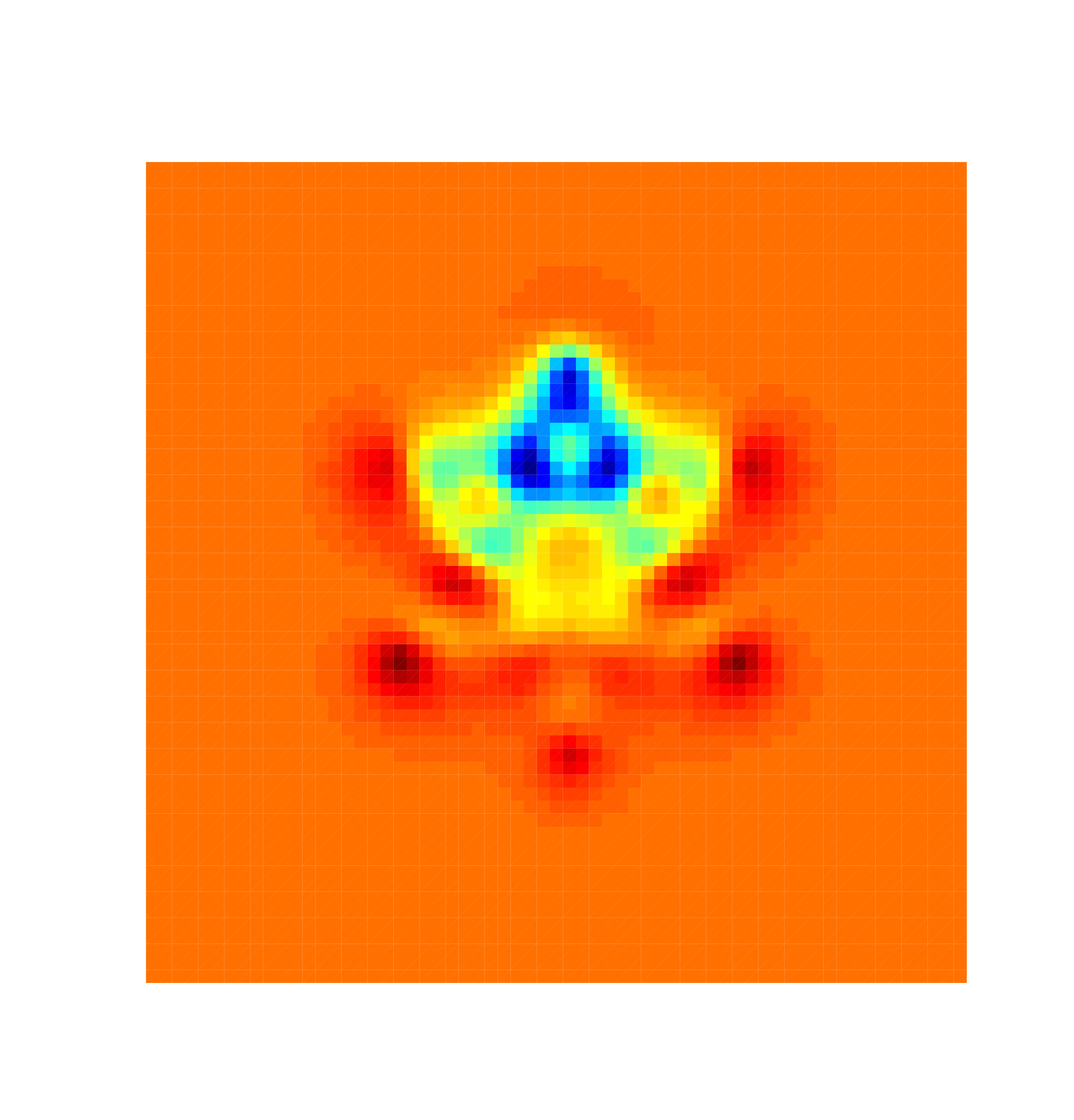} 
    \caption{Slice of $\hat{\chi}_{0}(0)[\hat{g}](\vr).$}
    \label{fig:slice}  
  \end{subfigure}
	\caption{Slices $\hat{g}$ and $\hat{\chi}_{0}(0)[\hat{g}](\vr)$ at the vertical midpoint of domain}
\end{figure}

We constructed $\hat{g}$ as a Gaussian centered at one of the carbon
atoms in the benzene molecule. Fig.~\ref{fig:g} shows a slice along
the $z$ direction of the
function $\hat{g},$ and~\ref{fig:V} shows the corresponding slice of potential function
$V$ for the benzene
molecule.  The pole expansion for computing the required quantities in
\eqref{eqn:multishiftOriginal} took 4400 seconds. This time encompasses
solving 15 sets of systems each with 200 shifts of the form
$\varepsilon_{i}+i\eta_{l}$, and each $\varepsilon_{i}$ corresponds to a distinct right
hand side. Computing the 15 smallest eigenvectors took 46.86 seconds.
Conversely, to use the alternative formula in
\eqref{eqn:multishiftOriginal} for computing
$\hat{\chi}_{0}(i\omega)[\hat{g}](\vr)$ requires computing a large
number of additional eigenvectors. Even just computing the first 1000
eigenvectors took 6563 seconds using LOBPCG. Computing the first
2000 eigenvectors took 82448 seconds. Fig.~\ref{fig:slice} shows an
example of the solution in a slice of the domain for $\omega = 0.$ 

\subsection{Non-orthogonal basis functions: $S\ne I$}
\begin{figure}
  \centering
  \begin{subfigure}[b]{.49\textwidth}
  	\begin{center}
	 	\includegraphics[width=.7\textwidth]{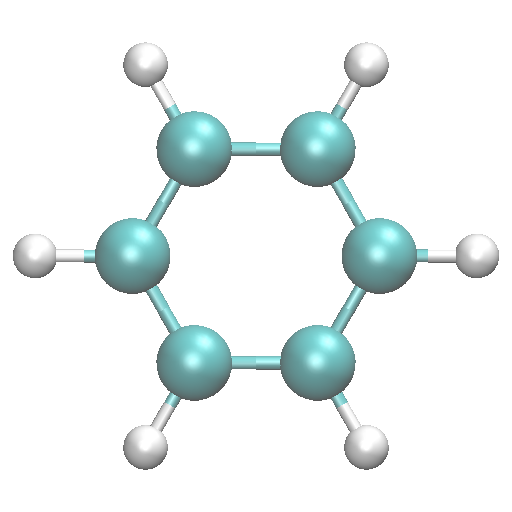}
	 	\caption{}
	 	\label{fig:benzene}
 	\end{center}
  \end{subfigure}
  \begin{subfigure}[b]{.49\textwidth}
  	\begin{center}
	  	\includegraphics[width=1\textwidth]{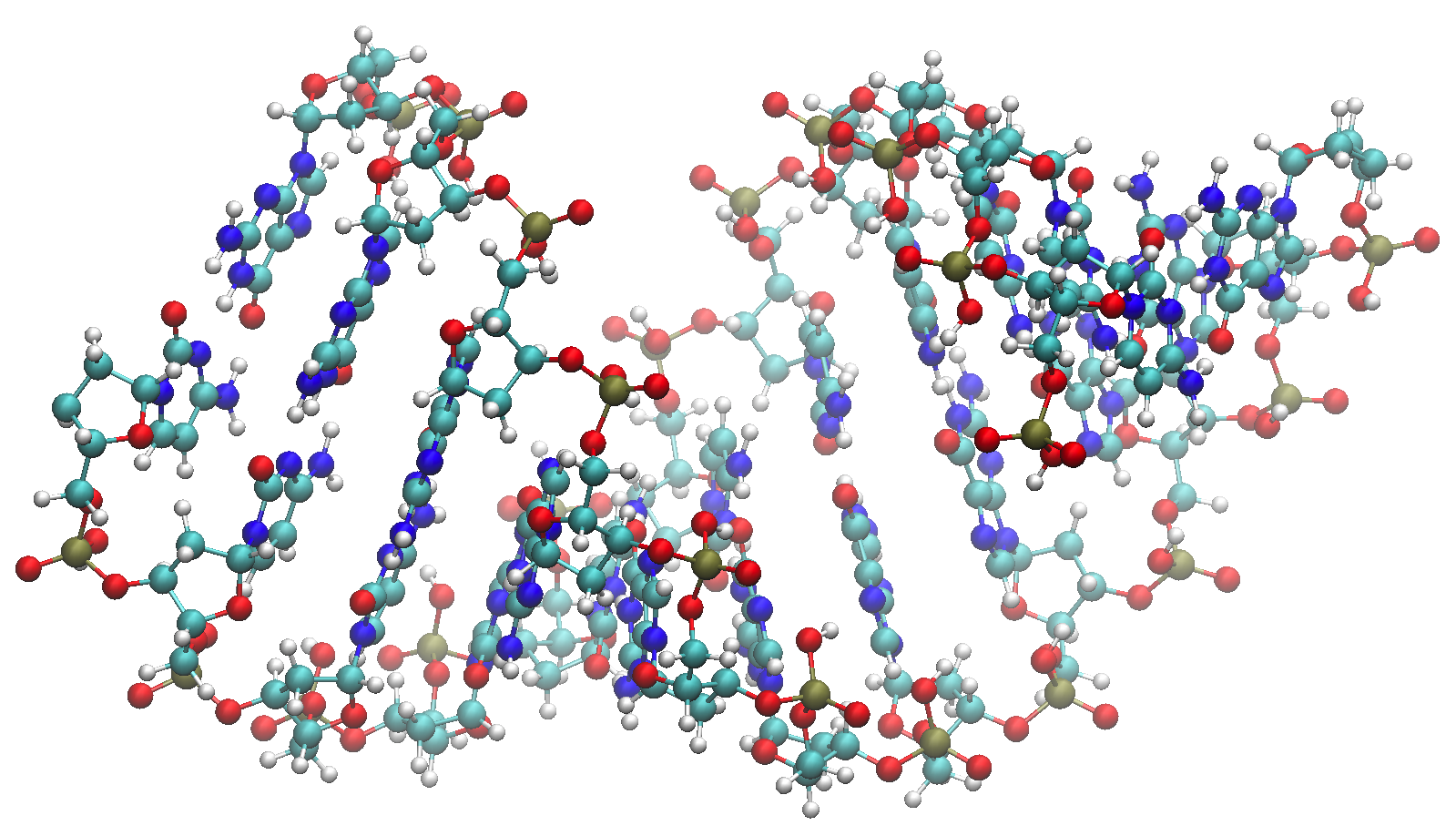}
	  	\caption{}
	  	\label{fig:dna}
  	\end{center}
  \end{subfigure}
   \caption{(a) The atomic configuration of a benzene molecule with $12$ atoms. (b) The atomic configuration of a DNA molecule with $715$ atoms.}
\end{figure}

We now consider the case where non-orthogonal basis functions are used
in the formulation of the problems as described in Section~\ref{sec:connection}. The matrices $H$ and $S$ are obtained by solving the Kohn-Sham density functional
theory problem for a DNA molecule with $715$ atoms using the SIESTA
package~\cite{SolerArtachoGaleEtAl2002} using the atomic orbital basis. 
The atomic configuration of the DNA molecule is given in Fig.~\ref{fig:dna}.
The number of electrons is $2442$ and the number of
occupied states $N_{e}=1221$ (spin degeneracy of $2$ is counted here). The $1221$ negative eigenvalues and the corresponding eigenvectors were directly computed and were only used to ensure that the right hand sides for the set of equations we solve are orthogonal to the negative eigenspace of the $(H,S)$ \response{pencil}.

In this case the basis functions are non-orthogonal and therefore $S\ne I$. This means that we are now interested in solving problems of the form
\begin{equation}
\label{eqn:exeq}
(H - i\eta_l S)x_l = b,
\end{equation} 
for a large number of equispaced $\eta_l$ in the interval $[-10,10].$ Furthermore, we are interested in solving the system for the same set of $\eta_l$ for multiple right hand sides $b.$ In this case both $H$ and $S$ are sparse and of size $7752 \times 7752.$ Therefore, the problem lends itself to the use of a direct method when solving the sub-problems for the pole expansion method as discussed in~\ref{sec:pole}. Since once again the problem is real, this means that we have to compute $P/2$ LU factorizations of matrices of the form
\begin{equation}
\label{eqn:exsub}
(H - \xi_k S).
\end{equation}
After this step has been completed the pole expansion method may be used to quickly solve problems of the form~\eqref{eqn:exeq}. First, we may compute the set of weights for each $\eta_l.$ Then, we may use the LU factorizations of~\eqref{eqn:exsub} to find the vectors $\tilde{h}_k.$ Finally, we may compute the solutions to the set of problems~\eqref{eqn:exeq} for a fixed $b$ using the computed weights. For each additional right had side we just need to compute a new set of $\tilde{h}_k$ and then combine them using the same weights as before. The marginal cost for each additional right had side is $P /2$ forward and backward substitutions plus the cost of combining $P$ vectors with $N_z$ distinct sets of weights. Comparatively, using the Lanczos method to solve for multiple right hand sides requires starting over because the Krylov subspace is dependent on $b.$ It is important to note that if the number of shifts is less than the number of poles used it would be more efficient to just factor the $N_z$ shifted systems. However, we are interested in the case where there are more shifts than poles even though some cases in the example do not reflect this situation.

This solution strategy is not dependent on computing LU factorizations
of~\eqref{eqn:exsub}. Any direct method may be used that allows for
rapid computation of solutions given a new right hand side. Therefore,
the pole expansion algorithm used in this manner has two distinct parts.
There is the factorization step where $P$ factorizations are computed,
and there is the solving step, where the factorizations are used to
solve $P$ sub-problems for each right hand side and then the solutions
are combined for all the desired $\eta_l.$ For a small number of shifts
and for very few right hand sides we still expect the Lanczos method to
potentially outperform the pole expansion method. However, as soon as we
have a large number of shifts, or there are enough right hand sides to
make the use of the direct methods favorable to an iterative method we
expect the pole expansion method to take much less time than the Lanczos
method.  

Based on the splitting of the work for the pole expansion between a
factorization step and a solve step, we present the results for this
example slightly differently than before. The problems were solved for a
varying number of $\eta_l,$ denoted $N_z.$ For each set of $\eta_l$ the
problem was solved using both the pole expansion method and the Lanczos
method for $N_{\text{rhs}}$ distinct right hand sides via the
transformation in~\eqref{eqn:chol_eq} and we report the average time to
solve the problem for a single right hand side along with the time for
computing the Cholesky factorization of $S,$ denoted $T_c.$ Table
\ref{tab:lanczos} shows the results of using the Lanczos method for
solving the set of problems. We observe that the number of iterations is
consistent regardless of the number of shifts, and because the problem
is small the method scales well as the number of shifts grows. However,
for each right hand side the Lanczos method must start from scratch.

Table~\ref{tab:pe} shows the time taken to factor the $P /
2$ sub-problems of the pole expansion method and then the cost for
computing a solution for all the $\eta_l$ per right hand side. Similar to before, we did not take advantage of the parallelism in the method and simply computed the factorizations sequentially. Here we
observe that the bulk of the computation time is in the factorization
step, which is expected since this now behaves like a direct method.
However, once the factorizations have been computed the marginal cost of
forming the solutions for all of the desired shifts with a new right
hand side is minimal. 

\begin{table}
  \centering
  \begin{tabular}{c c | c | c | c | c |  }
    \cline{3-6}
    & & \multicolumn{4}{|c|}{Lanczos} \\ \hline
    \multicolumn{1}{|c|}{ $N_z$ } &  \multicolumn{1}{|c|}{ $N_{\text{rhs}}$  } & Avg. Iter. & $T_c(s)$ & Avg. Solve Time(s) per $b$ & Max Error  \\ \hline
    \multicolumn{1}{|c|}{3} & \multicolumn{1}{|c|}{20} & 171 & 1.33  & 102  & 3.88 $\times 10^{-9}$   \\ \hline
    \multicolumn{1}{|c|}{11} & \multicolumn{1}{|c|}{20} & 173 & 1.33  & 103.8  & 4.03 $\times 10^{-9}$   \\ \hline
    \multicolumn{1}{|c|}{101} & \multicolumn{1}{|c|}{20} & 173 & 1.33  & 114.8  & 4.25 $\times 10^{-9}$   \\ \hline
    \multicolumn{1}{|c|}{1001} & \multicolumn{1}{|c|}{20} & 171 & 1.33  & 223.6  & 4.94 $\times 10^{-9}$   \\ \hline
  \end{tabular}
  \caption{Lanczos method for solving $(H-i\eta_lS)x_l = b.$}
  \label{tab:lanczos}
\end{table}

\begin{table}
  \centering
  \begin{tabular}{c c  | c | c | c | c |  }
    \cline{3-6}
    & & \multicolumn{4}{|c|}{Pole Expansion} \\ \hline
     \multicolumn{1}{|c|}{ $N_z$ } &\multicolumn{1}{|c|}{ $N_{\text{rhs}}$ } &  P & Factor Time(s) & Avg. Solve Time(s) per $b$ & Max Error \\ \hline
    \multicolumn{1}{|c|}{3} & \multicolumn{1}{|c|}{20} & 60 & 149.4 & 3.95 & 5.52 $\times 10^{-10}$  \\ \hline
    \multicolumn{1}{|c|}{11}  & \multicolumn{1}{|c|}{20} & 60 & 149.4 & 3.98  & 5.46 $\times 10^{-10}$ \\ \hline
    \multicolumn{1}{|c|}{101}  & \multicolumn{1}{|c|}{20} & 60 & 149.4 & 5.14  & 5.45 $\times 10^{-10}$ \\ \hline
    \multicolumn{1}{|c|}{1001}  & \multicolumn{1}{|c|}{20} & 60 & 149.4 & 16.2  & 5.56 $\times 10^{-10}$ \\ \hline
  \end{tabular}
  \caption{Pole expansion method for solving $(H-i\eta_lS)x_l = b.$}
  \label{tab:pe}
\end{table}

When a direct method is an option for solving problems of the form
\eqref{eqn:exeq} and there are a large number of shifts, the direct
method may be combined with the pole expansion method to essentially
parametrize the factorizations of $N_z$ distinct matrices on the
factorizations of $P$ distinct matrices. If the number of shifts is much
larger than the number of poles required for the desired accuracy, this
reduction in the number of required factorizations turns out to be very
beneficial computationally. Furthermore, if memory is an issue, it is
possible to only ever store one factorization at a time. Specifically,
once a factorization is computed for a given pole, the vectors
$\tilde{h}_k$ may be computed for each right hand side and then the
factorization may be discarded. Overall, the combined use of the pole
expansion and an efficient direct method appears to be a very efficient
method for solving sets of parametrized linear systems for multiple
right hand sides, or, in some cases, even for a single right hand side.


\section{Conclusion}\label{sec:conclusion}
We have presented a new method for efficiently solving a type of
parametrized linear systems, motivated from electronic structure
calculations. By building a quadrature scheme based on the ideas in
\cite{HHT} we are able to represent the solutions of the parametrized
shifted systems as weighted linear combinations of solutions to a set of
fixed problems, where the weights vary based on the parameter of the
system. This method scales well as the number of distinct parameters for
which we want to solve~\eqref{eqn:multishift} grows. Furthermore,
because the solutions to the parametrized equations are based on a fixed
set of sub-problems there is flexibility in how the sub-problems are
solved. We presented examples using both iterative and direct solvers
within the framework for solving the shifted systems. The method
presented here can be more favorable compared to a Lanczos based method,
especially when solutions to a large number of parameters or a large
number of right hand sides are required. 



\section*{Acknowledgments}

A. D. is currently supported by NSF Fellowship DGE-1147470 and was partially supported by the National Science Foundation grant DMS-0846501.
L. L. was partially supported by the Laboratory Directed Research and
Development Program of Lawrence Berkeley National Laboratory under the
U.S.  Department of Energy contract number DE-AC02-05CH11231, and by
Scientific Discovery through Advanced Computing (SciDAC) program funded
by U.S.  Department of Energy, Office of Science, Advanced Scientific
Computing Research and Basic Energy Sciences. L. Y. was partially
supported by National Science Foundation under award DMS-0846501 and by
the Mathematical Multifaceted Integrated Capability Centers (MMICCs)
effort within the Applied Mathematics activity of the U.S. Department of
Energy's Advanced Scientific Computing Research program, under Award
Number(s) DE-SC0009409. We would like to thank Lenya Ryzhik for providing computing resources.  We are grateful to Alberto Garcia and Georg Huhs for providing the atomic configuration for the DNA molecule.  L. L. thanks Dario Rocca and Chao Yang for helpful discussions, and thanks the
hospitality of Stanford University where the idea of this paper started.


\bibliographystyle{siam}
\bibliography{multishift}

\begin{thebibliography}{10}

\bibitem{Adler1962}
{\sc S.~L. Adler}, {\em Quantum theory of the dielectric constant in real
  solids}, Phys. Rev., 126 (1962), pp.~413--420.

\bibitem{AryasetiawanGunnarsson1998}
{\sc F.~Aryasetiawan and O.~Gunnarsson}, {\em The {GW} method}, Rep. Prog.
  Phys., 61 (1998), p.~237.

\bibitem{BaiFreund2001}
{\sc Z.~Bai and Roland~W. Freund}, {\em A partial {P}ad{\'e}-via-{L}anczos
  method for reduced-order modeling}, Linear Algebra Appl., 332 (2001),
  pp.~139--164.

\bibitem{BaroniGironcoliDalEtAl2001}
{\sc S.~Baroni, S.~de~Gironcoli, A.~Dal~Corso, and P.~Giannozzi}, {\em Phonons
  and related crystal properties from density-functional perturbation theory},
  Rev. Mod. Phys., 73 (2001), pp.~515--562.

\bibitem{BaroniGiannozziTesta1987}
{\sc S.~Baroni, P.~Giannozzi, and A.~Testa}, {\em Green's-function approach to
  linear response in solids}, Phys. Rev. Lett., 58 (1987), pp.~1861--1864.

\bibitem{DattaSaad1991}
{\sc B.~N. Datta and Y.~Saad}, {\em {Arnoldi} methods for large
  {Sylvester-like} observer matrix equations, and an associated algorithm for
  partial spectrum assignment}, Linear Algebra Appl., 154 (1991), pp.~225--244.

\bibitem{DriscollSC1}
{\sc T.~A. Driscoll}, {\em Algorithm 756: a {MATLAB} toolbox for
  {S}chwarz-{C}hristoffel mapping}, ACM Trans. Math. Software, 22 (1996),
  pp.~168--186.

\bibitem{DriscollSC2}
\leavevmode\vrule height 2pt depth -1.6pt width 23pt, {\em Algorithm 843:
  Improvements to the {S}chwarz-{C}hristoffel toolbox for {MATLAB}}, ACM Trans.
  Math. Software, 31 (2005), pp.~239--251.

\bibitem{FeldmannFreund1995}
{\sc P.~Feldmann and R.~W. Freund}, {\em Efficient linear circuit analysis by
  {P}ad{\'e} approximation via the {L}anczos process}, IEEE Trans. Comput. Aid
  D., 14 (1995), pp.~639--649.

\bibitem{TFQMR}
{\sc R.~Freund}, {\em A transpose-free quasi-minimal residual algorithm for
  non-{H}ermitian linear systems}, SIAM J. Sci. Comput., 14 (1993),
  pp.~470--482.

\bibitem{FriedrichSchindlmayr2006}
{\sc C.~Friedrich and A.~Schindlmayr}, {\em Many-body perturbation theory: the
  {GW} approximation}, NIC Series, 31 (2006), p.~335.

\bibitem{Frommer2003}
{\sc A.~Frommer}, {\em {BiCGStab(l)} for families of shifted linear systems},
  Computing, 70 (2003), pp.~87--109.

\bibitem{FrommerGlassner1998}
{\sc A.~Frommer and U.~Gl{\"a}ssner}, {\em Restarted {GMRES} for shifted linear
  systems}, SIAM J. Sci. Comput., 19 (1998), pp.~15--26.

\bibitem{Furche2001}
{\sc F.~Furche}, {\em Molecular tests of the random phase approximation to the
  exchange-correlation energy functional}, Phys. Rev. B, 64 (2001), p.~195120.

\bibitem{GallivanGrimmeVan1996}
{\sc K.~Gallivan, G.~Grimme, and P.~Van~Dooren}, {\em A rational {L}anczos
  algorithm for model reduction}, Numer. Algorithms, 12 (1996), pp.~33--63.

\bibitem{GiustinoCohenLouie2010}
{\sc F.~Giustino, M.~L. Cohen, and S.~G. Louie}, {\em {GW} method with the
  self-consistent sternheimer equation}, Phys. Rev. B, 81 (2010), p.~115105.

\bibitem{GVL}
{\sc G.~H. Golub and C.~F. Van~Loan}, {\em Matrix computations}, Johns Hopkins
  Univ. Press, Baltimore, third~ed., 1996.

\bibitem{GonzeAllanTeter1992}
{\sc X.~Gonze, D.~C Allan, and M.~P. Teter}, {\em Dielectric tensor, effective
  charges, and phonons in $\alpha$-quartz by variational density-functional
  perturbation theory}, Phys. Rev. Lett., 68 (1992), p.~3603.

\bibitem{HHT}
{\sc N.~Hale, N.~J. Higham, and L.~N. Trefethen}, {\em Computing
  \protect{$A^{\alpha}$}, \protect{$\log(A)$}, and related matrix functions by
  contour integrals}, SIAM J. Numer. Anal., 46 (2008), pp.~2505--2523.

\bibitem{Hedin1965}
{\sc L.~Hedin}, {\em New method for calculating the one-particle {Green's}
  function with application to the electron-gas problem}, Phys. Rev., 139
  (1965), p.~A796.

\bibitem{HohenbergKohn1964}
{\sc P.~Hohenberg and W.~Kohn}, {\em {Inhomogeneous electron gas}}, Phys. Rev.,
  136 (1964), pp.~B864--B871.

\bibitem{Knyazev2001}
{\sc A.~V. Knyazev}, {\em Toward the optimal preconditioned eigensolver:
  Locally optimal block preconditioned conjugate gradient method}, SIAM J. Sci.
  Comp., 23 (2001), p.~517.

\bibitem{KohnSham1965}
{\sc W.~Kohn and L.~Sham}, {\em {Self-consistent equations including exchange
  and correlation effects}}, Phys. Rev., 140 (1965), pp.~A1133--A1138.

\bibitem{Lanczos1950}
{\sc C.~Lanczos}, {\em An iteration method for the solution of the eigenvalue
  problem of linear differential and integral operators}, J. Res. Nat. Bur.
  Stand., 45 (1950), pp.~255--282.

\bibitem{LangrethPerdew1975}
{\sc D.~C. Langreth and J.~P. Perdew}, {\em The exchange-correlation energy of
  a metallic surface}, Solid State Commun., 17 (1975), pp.~1425--1429.

\bibitem{LangrethPerdew1977}
\leavevmode\vrule height 2pt depth -1.6pt width 23pt, {\em Exchange-correlation
  energy of a metallic surface: Wave-vector analysis}, Phys. Rev. B, 15 (1977),
  p.~2884.

\bibitem{LinLuYingE2009}
{\sc L.~Lin, J.~Lu, L.~Ying, and W.~E}, {\em Pole-based approximation of the
  {F}ermi-{D}irac function}, Chin. Ann. Math., 30B (2009), p.~729.

\bibitem{Martin2004}
{\sc R.~Martin}, {\em Electronic Structure -- Basic Theory and Practical
  Methods}, Cambridge Univ. Pr., West Nyack, {NY}, 2004.

\bibitem{MARROK05}
{\sc P.G. Martinsson and V.~Rokhlin}, {\em A fast direct solver for boundary
  integral equations in two dimensions}, J. Comput. Phys., 205 (2005).

\bibitem{Meerbergen2003}
{\sc K.~Meerbergen}, {\em The solution of parametrized symmetric linear
  systems}, SIAM J. Matrix Anal. Appl., 24 (2003), pp.~1038--1059.

\bibitem{NguyenGironcoli2009}
{\sc H.-V. Nguyen and S.~de~Gironcoli}, {\em Efficient calculation of exact
  exchange and {RPA} correlation energies in the adiabatic-connection
  fluctuation-dissipation theory}, Phys. Rev. B, 79 (2009), p.~205114.

\bibitem{PASA1975}
{\sc C.~Paige and M.~Saunders}, {\em Solution of sparse indefinite systems of
  linear equations}, SIAM Journal on Numerical Analysis, 12 (1975),
  pp.~617--629.

\bibitem{Parks2006}
{\sc M.~Parks, E.~de~Sturler, G.~Mackey, D.~Johnson, and S.~Maiti}, {\em
  Recycling krylov subspaces for sequences of linear systems}, SIAM Journal on
  Scientific Computing, 28 (2006), pp.~1651--1674.

\bibitem{PingRoccaGalli2013}
{\sc Y.~Ping, D.~Rocca, and G.~Galli}, {\em Electronic excitations in light
  absorbers for photoelectrochemical energy conversion: first principles
  calculations based on many body perturbation theory.}, Chem. Soc. Rev., 42
  (2013), pp.~2437--2469.

\bibitem{Saad1986}
{\sc Y.~Saad and M.~Schultz}, {\em Gmres: A generalized minimal residual
  algorithm for solving nonsymmetric linear systems}, SIAM Journal on
  Scientific and Statistical Computing, 7 (1986), pp.~856--869.

\bibitem{SABAKI2012}
{\sc A.~Saibaba, T.~Bakhos, and P.~Kitanidis}, {\em A flexible krylov solver
  for shifted systems with application to oscillatory hydraulic tomography},
  SIAM Journal on Scientific Computing, 35 (2013), pp.~A3001--A3023.

\bibitem{SimonciniSzyld2007}
{\sc V.~Simoncini and D.~B. Szyld}, {\em Recent computational developments in
  {Krylov} subspace methods for linear systems}, Numer. Linear Algebra Appl.,
  14 (2007), pp.~1--59.

\bibitem{SolerArtachoGaleEtAl2002}
{\sc J.~M. Soler, E.~Artacho, J.~D. Gale, A.~Garc{\'\i}a, J.~Junquera,
  P.~Ordej{\'o}n, and D.~S{\'a}nchez-Portal}, {\em {The SIESTA method for ab
  initio order-N materials simulation}}, J. Phys.: Condens. Matter, 14 (2002),
  pp.~2745--2779.

\bibitem{UmariStenuitBaroni2010}
{\sc P.~Umari, G.~Stenuit, and S.~Baroni}, {\em {GW} quasiparticle spectra from
  occupied states only}, Phys. Rev. B, 81 (2010), p.~115104.

\bibitem{Wiser1963}
{\sc N.~Wiser}, {\em Dielectric constant with local field effects included},
  Phys. Rev., 129 (1963), pp.~62--69.

\bibitem{YangMezaLeeEtAl2009}
{\sc C.~Yang, J.~C. Meza, B.~Lee, and L.~W. Wang}, {\em {KSSOLV--a MATLAB
  toolbox for solving the Kohn--Sham equations}}, ACM Trans. Math. Software, 36
  (2009), p.~10.

\end{thebibliography}

\end{document}